\documentclass[11pt, a4paper]{amsart}
\usepackage[utf8]{inputenc}
\usepackage[english]{babel}

\usepackage[margin=2.5cm]{geometry}
\usepackage{mathrsfs}
\usepackage{enumitem}
\usepackage{pinlabel}
\usepackage{graphicx}
\usepackage{amsmath,amsfonts,amssymb,mathtools}
\usepackage{tikz,tikz-cd}
\usepackage{subcaption}
\usepackage{hyperref}
\usepackage{soul}
\usetikzlibrary{calc,patterns,arrows,shapes.arrows,intersections}
\usetikzlibrary{decorations}

\newtheorem{thm}{Theorem}[section]

\newtheorem{prop}[thm]{Proposition}
\newtheorem{lem}[thm]{Lemma}
\newtheorem{cor}[thm]{Corollary}
\newtheorem{rem}[thm]{Remark}
\newtheorem{example}[thm]{Example} 
\newtheorem{q}[thm]{Question}

\newtheorem{defn}[thm]{Definition}

\newcommand{\bN}{\mathbb{N}}
\newcommand{\bZ}{\mathbb{Z}}
\newcommand{\bK}{\mathbb{K}}
\newcommand{\bH}{\mathbb{H}}
\newcommand{\bF}{\mathbb{F}}
\newcommand{\cF}{\mathcal{F}}
\newcommand{\cH}{\mathcal{H}}

\newcommand{\tG}{{\tt G}}
\newcommand{\tH}{{\tt H}}
\newcommand{\tC}{{\tt C}}
\newcommand{\tP}{{\tt P}}
\newcommand{\tK}{{\tt K}}
\newcommand{\tL}{{\tt L}}
\newcommand{\tW}{{\tt W}}

\newcommand{\tT}{{\tt T}}
\newcommand{\Hh}{\mathrm{H}^h}
\newcommand{\UC}{\ddot{\mathrm{C}}}
\newcommand{\UH}{\ddot{\mathrm{H}}}
\newcommand{\cl}{{\rm Col}}
\newcommand{\e}{\varepsilon}

\definecolor{aquamarine}{rgb}{0.5, 1.0, 0.83}
\definecolor{princetonorange}{rgb}{1.0, 0.56, 0.0}
\definecolor{caribbeangreen}{rgb}{0.0, 0.8, 0.6}
\definecolor{bunired}{rgb}{0.8, 0.0, 0.0}
\definecolor{cdgreen}{rgb}{0.0, 0.42, 0.24}
\definecolor{lavender(floral)}{rgb}{0.71, 0.49, 0.86}

\title{Categorifying connected domination via graph \"uberhomology}

\author{Luigi Caputi}
\author{Daniele Celoria}
\author{Carlo Collari}
\date{}

\begin{document}

\begin{abstract}
\"Uberhomology is a recently defined homology theory for simplicial complexes, which yields subtle information on graphs.
We prove that bold homology, a certain specialisation of \"uberhomology, is related to dominating sets in graphs. To this end, we interpret \"uberhomology as a poset homology, and investigate its functoriality properties. 
We then show that the Euler characteristic of the bold homology of a graph coincides with an evaluation of its connected domination polynomial. Even more, the bold chain complex retracts onto a complex generated by connected dominating sets.
We conclude with several computations of this homology on families of graphs; these include a vanishing result for trees, and a characterisation result for complete graphs. 
\end{abstract}

\maketitle

\section{Introduction}

The purpose of this paper is twofold; first, we investigate the functoriality properties of the recently defined  \"uberhomology \cite{uberhomology}. We then provide a categorification of an evaluation of the connected domination polynomial~$D^c_{\tG}(x)$.
Unexpectedly, these two directions turn out to be closely related; in fact, we show that  $D^c_{\tG}(-1)$ is the Euler characteristic of a suitable homology theory~$\bH(\tG)$. This latter homology is a degree specialisation of the \"uberhomology, but it also admits an independent definition, \emph{cf.}~ \cite[Section 8]{uberhomology}.

Recall that the \"uberhomology $\UH(X)$ is a combinatorially defined (triply-graded) homology theory associated to a finite and connected simplicial complex~$X$.
Its definition relies on certain combinatorial filtrations on the simplicial chain complex of $X$, arranged in a poset-like fashion reminiscent of Khovanov homology~\cite{Khovanov}. This is not a coincidence; indeed, we show that the \"uberhomology is a special case of a poset homology (in the sense of \cite{chandler2019posets,primo}, \emph{cf}.~Remark~\ref{rmk:posethom}).
A consequence of this interpretation yields the following result:

\begin{thm}\label{thm:functor uberhom}
The \"uberhomology is a bi-functor
\[
\UH^*(-;-)\colon \mathbf{RegSimpl} \times {\bf Ring} \to \mathbf{grAb}
\]
where $\mathbf{RegSimpl}$ denotes the category of simplicial complexes and injective simplicial maps.
\end{thm}

\"Uberhomology groups measure both combinatorial and topological features of simplicial complexes and therefore, in particular, of simple graphs. 
We focus our attention on the \"uberhomology in a specific bi-degree,
namely we define the \emph{bold homology}~$\bH_*(X) $ of $X$  as  $\UH_{0,0}^*(X)$. 
For a graph $\tG$, the homology $\bH(\tG)$ has an independent definition as the homology of the chain complex ${\rm C}\bH(\tG)$ -- \emph{cf.}~\cite[Secion~8]{uberhomology}. 
A basis of $\bH(\tG)$ is provided by subgraphs of $\tG$.
We use these facts, in conjunction with Theorem~\ref{thm:functor uberhom}, to prove our main result:
\begin{thm}\label{thm:categorification}
The bold homology $\bH(\tG)$ is a categorification of $D^c_{\tG}(-1)$; that is $\bH$ is functorial with respect to injective morphisms of graphs, and its Euler characteristic is~$D^c_{\tG}(-1)$.
\end{thm}

Recall that, for a graph $\tG$, a dominating set $D$ is a subset of the vertices of $\tG$ such that every vertex in $\tG$ is in $D$ or adjacent to at least one member of $D$. 
Dominating sets in graphs have been extensively studied (see~\emph{e.g.}~\cite{haynes2013fundamentals} and references therein); finding dominating sets of a given size is well known to be a NP-complete problem, and is related to open conjectures (such as Vizing's conjecture~\cite{Vizing_1968}, see also the survey \cite{vizingconj}).

The categorification provided by Theorem~\ref{thm:categorification} can be strengthened to uncover a deeper relation between connected dominating sets and $\bH$.
More precisely, if ${\rm C}\bH(\tG)$  denotes the chain complex computing the bold homology $\bH(\tG)$, we obtain the following result:

\begin{thm}
There exists a quasi-isomorphism between the chain complex ${\rm C}\bH(\tG)$ and a complex~${\rm D}\bH(\tG)$. The complex ${\rm D}\bH(\tG)$ is spanned by connected dominating sets of $\tG$, and its differential is induced by inclusions. 
\end{thm}

The proof of this last theorem relies on various techniques from combinatorial algebraic topology (see~\cite{Kozlov}), and, in particular, from algebraic Morse theory. We prove some technical results on discrete Morse matchings (\emph{cf}.~Lemma~\ref{lem:tech lemma} and Remark~\ref{rmk:induced diff koz}), which might be of independent interest.
These techniques enable us to compute the bold homology for certain families of graphs; namely, trees, complete bipartite graphs, and cycle graphs, proving~\cite[Conjecture~8.2]{uberhomology}. We also examine the behaviour of $\bH$ under certain natural graph operations -- \emph{cf.}~Propositions~\ref{prop:cone} and~\ref{prop:neckstretch}.
Finally, we prove that the bold homology characterises complete graphs:

\begin{thm}
The homology $\bH_1(\tG)$ is non-zero if and only if $\tG$ is a complete graph.
\end{thm}

We point out that computations of~$\bH(\tG)$ can be carried out by means of computer software -- see the Sage~\cite{sagemath}  implementation~\cite{githububer}.\\
We conclude with some sample computations of sporadic examples (see Table~\ref{tab:summary_table}), and a list of open questions.

\subsection*{Organisation of the paper}
In the first section, we recall the definition of \"uberhomology, and provide an alternative interpretation using poset homology with functor coefficients. This viewpoint is put to use in Section~\ref{sec:functoriality}, where we investigate the functoriality of the \"uberhomology with respect to injective morphisms of graphs. In Section~\ref{sec: bold homology and dom} we recall the definition of the bold homology, as a specialisation of  \"uberhomology, and prove that its Euler characteristic categorifies an evaluation of the connected domination polynomial. We provide some applications and computations in Section~\ref{sec: appl and comps}, and conclude with some open questions.

\subsection*{Conventions}
Typewriter font, e.g.~$\tG$, $\tH$, \emph{etc.}, are used to denote finite, simple and connected graphs, possibly oriented. 
Unless otherwise stated, $R$ denotes a ring, $\bK$ is a field,  and $\bF$ denotes the field with two elements. 
We include here the notation for some graph families used in the paper:
$\tL_n$ denotes the linear graph,
 $\tC_n$ the  cycle,
 $\tW_{n}$ the  wheel graph, and $\tK_{n}$ the complete graph on $n$ vertices. We also denote by  $\tK_{m,n}$ the complete bipartite graph on $m, n$ vertices, and by
${\tt Cube}(n)$ the $1$-skeleton of the $n$-dimensional cube.

\subsection*{Acknowledgements}

The authors are thankful to F.~Petrov for suggesting the proof of Proposition~\ref{prop:russo}. LC~acknowledges support from the \'{E}cole Polytechnique F\'{e}d\'{e}rale de Lausanne via a collaboration agreement with the University of Aberdeen.
DC was partially supported by the European Research Council (ERC) under the EU Horizon 2020 research and innovation programme (grant agreement No 674978) and by Hodgson-Rubinstein's ARC grant DP190102363 ``Classical And Quantum Invariants Of Low-Dimensional Manifolds''.
During the writing of this paper CC was a postdoc at the New York University Abu Dhabi.

\section{\"Uber and poset homologies}

In this section we recall some basic notions and prove that the \"uberhomology is a poset homology.

\subsection{\"Uberhomology}\label{sec:uberhom}

We start by giving a brief account of the definitions from~\cite{uberhomology}. 

Let $X$ be a finite and connected simplicial complex with $m$ vertices, which we assume to be ordered, say $V(X)=\{v_1,\dots,v_m\}$. 

\begin{defn}\label{df:bicol simpl compl}
A \emph{bi-colouring} $\varepsilon$ on $X$ is a map $\varepsilon\colon V(X) \to \{0,1\}$. A \emph{bi-coloured simplicial complex } is a pair $(X,\e)$ consisting of a simplicial complex $X$ and a {bi-colouring} $\varepsilon$ on $V(X)$.
\end{defn}

Given a $n$-dimensional simplex~$\sigma$ in a bi-coloured simplicial complex $(X,\e)$, define its \emph{weight} with respect to $\e$ as the sum
\begin{equation}\label{eq:weight}
w_\varepsilon(\sigma )\coloneqq \dim(\sigma) +1-\sum_{v_i\in V(\sigma)} \varepsilon(v_i) \ .
\end{equation}
In other words, the weight is the number of $0$-coloured vertices in a simplex.
If we fix a colouring~$\e$, the weight in Equation~\eqref{eq:weight} induces a filtration of the simplicial chain complex~$C_*(X;R)$ associated to~$X$. More explicitly, we set
\[ \mathscr{F}_j (X, \e ;R)  \coloneqq R \langle\ \sigma \mid w_{\e}(\sigma) \leq j\ \rangle \subseteq C_*(X;R) \ . \]
The simplicial differential $\partial$ preserves this filtration, so each~$ \mathscr{F}_j (X, \e ;R) $ is a sub-complex of~$C_*(X;R)$. We can decompose $\partial$ as the sum of two differentials (\emph{cf}.~\cite[Lemma~2.1]{uberhomology}); one which preserves the weight,  and one which decreases it by one.  We denote the former by $\partial_h$, and won't make use of the latter. Call $(C(X,\e;R),\partial_h)$ the bi-graded chain complex whose underlying module  is $C(X;R)$; the first degree is given by simplices' dimensions, and the second is given by the weight $w_\e$.

\begin{defn}\label{def:hor hom}
The $\e$-\emph{horizontal homology} $\Hh(X,\e;R)$ of $(X,\e)$ is the homology of the bi-graded chain complex $(C(X,\e;R),\partial_h)$.
\end{defn}

When the ring of coefficients $R$ is clear from the context, we will simply denote the $\e$-horizontal homology of $(X,\e)$ by $\Hh(X,\e)$. 

Consider the Boolean poset $B(m)$ on $m$ vertices, that is the set consisting of subsets of $\{1,\dots,m\}$, partially ordered by inclusion. 

\begin{rem}\label{rem:boolean col}
One can decorate the elements of~$B(m)$ with the bi-colourings on $X$. Indeed, the set of bi-colourings on $X$ can be canonically identified with the elements of $\{ 0 , 1\}^m$ via the map~$\e \mapsto (\e(v_1),...,\e(v_m))$. The minimum of $B(m)$ corresponds to the $(0, \ldots,0)$-colouring, and the maximum of $B(m)$ corresponding to the $(1,\ldots,1)$-colouring.
\end{rem}

We are now ready to recall the definition of the \"uberhomology -- \emph{cf}.~\cite[Section~6]{uberhomology}.
Let $\e$ and $\e'$ be two bi-colourings on $X$ which differ only on a vertex $v_i$; assume further that $\e(v_i) = 0$ and~$\e'(v_i) =1$.  We denote by $d_{\e, \e'} $ the weight-preserving part of the identity map $\mathrm{Id} \colon \mathrm{H}^h(X,\e) \to \mathrm{H}^h(X,\e')$. More explicitly
\[ d_{\e, \e'} (\sigma) = \begin{cases} \sigma & \text{if }w_{\e}(\sigma) = w_{\e'}(\sigma) \\ 0 &\text{otherwise.}\end{cases}  \]
Note that the second case can only occur if $w_{\e}(\sigma) = w_{\e'}(\sigma)-1$.

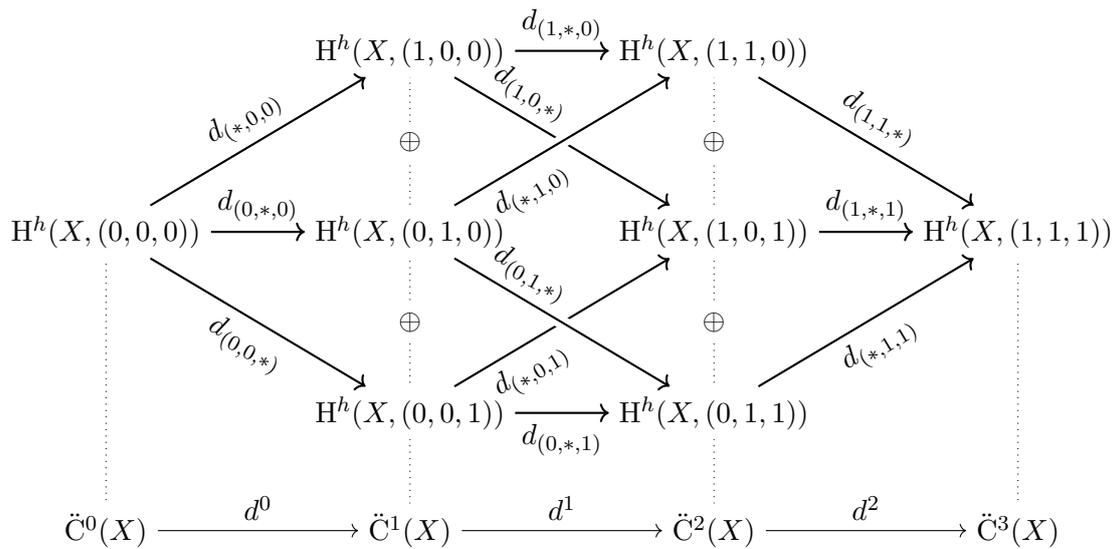
\begin{figure}[ht]

\begin{tikzpicture}[scale = .8]

\draw[dotted] (0,0) -- (0,-4.5) node[below] {$\UC^0(X)$};

\draw[dotted] (5,3) -- (5,-4.5) node[below] {$\UC^1(X)$};

\draw[dotted] (10,3) -- (10,-4.5) node[below] {$\UC^2(X)$};

\draw[dotted] (15,0) -- (15,-4.5) node[below] {$\UC^3(X)$};

\node[below] (uc0) at (0,-4.5) {\phantom{$\UC^0(X)$}};
\node[below] (uc1) at (5,-4.5) {\phantom{$\UC^0(X)$}};
\node[below] (uc2) at (10,-4.5) {\phantom{$\UC^0(X)$}};
\node[below] (uc3) at (15,-4.5) {\phantom{$\UC^0(X)$}};

\draw[->] (uc0) -- (uc1) node[midway, above] {$d^0$};
\draw[->] (uc1) -- (uc2) node[midway, above] {$d^1$};
\draw[->] (uc2) -- (uc3) node[midway, above] {$d^2$};

\node[fill, white] at (5,1.5){$\oplus$};
\node[fill, white] at (5,-1.5){$\oplus$};
\node  at (5,1.5){$\oplus$};
\node  at (5,-1.5){$\oplus$};

\node[fill, white] at (10,1.5){$\oplus$};
\node[fill, white] at (10,-1.5){$\oplus$};
\node  at (10,1.5){$\oplus$};
\node  at (10,-1.5){$\oplus$};

\node[fill, white] at (0,0) {${\Hh(X,(0,0,0))}$};

\node[fill, white] at (5,3){${\Hh(X,(1,0,0))}$};
\node[fill, white] at (5,0){${\Hh(X,(0,1,0))}$};
\node[fill, white] at (5,-3){${\Hh(X,(0,0,1))}$};

\node[fill, white] at (10,3) {${\Hh(X,(1,1,0))}$};
\node[fill, white] at (10,0) {${\Hh(X,(1,0,1))}$};
\node[fill, white] at (10,-3) {${\Hh(X,(0,1,1))}$};

\node[fill, white] at (15,0) {${\Hh(X,(1,1,1))}$};

\node (a) at (0,0) {${\Hh(X,(0,0,0))}$};

\node (b1) at (5,3) {${\Hh(X,(1,0,0))}$};
\node (b2) at (5,0) {${\Hh(X,(0,1,0))}$};
\node (b3) at (5,-3){${\Hh(X,(0,0,1))}$};

\node (c1) at (10,3) {${\Hh(X,(1,1,0))}$};
\node (c2) at (10,0) {${\Hh(X,(1,0,1))}$};
\node (c3) at (10,-3) {${\Hh(X,(0,1,1))}$};

\node (d) at (15,0) {${\Hh(X,(1,1,1))}$};

\draw[thick, ->] (a) -- (b1) node[midway,above,rotate =31] {$d_{(*,0,0)}$}; 
\draw[thick, ->] (a) -- (b2) node[midway,above] {$d_{(0,*,0)}$}; 
\draw[thick, ->] (a) -- (b3) node[midway,below,rotate =-29] {$d_{(0,0,*)}$};

\draw[thick, ->] (b1) -- (c1) node[midway,above] {$d_{(1,*,0)}$}; 
\draw[thick, ->] (b1) -- (c2) node[midway,above left,rotate =-29] {$d_{(1,0,*)}$}; 

\draw[thick, ->] (b3) -- (c2) node[midway,below left,rotate =31] {$d_{(*,0,1)}$}; 
\draw[thick, ->] (b3) -- (c3) node[midway,below] {$d_{(0,*,1)}$}; 

\draw[line width = 5, white] (b2) -- (c1) ; 
\draw[line width = 5, white] (b2) -- (c3) ; 
\draw[thick, ->] (b2) -- (c1) node[midway,below left,rotate =31] {$d_{(*,1,0)}$} ; 
\draw[thick, ->] (b2) -- (c3) node[midway,above left,rotate =-29] {$d_{(0,1,*)}$}; 

\draw[thick, <-] (d) -- (c1) node[midway,above,rotate =-29] {$d_{(1,1,*)}$}; 
\draw[thick, <-] (d) -- (c2) node[midway,above] {$d_{(1,*,1)}$}; 
\draw[thick, <-] (d) -- (c3) node[midway,below,rotate =31] {$d_{(*,1,1)}$}; 
\end{tikzpicture}

\caption{The Boolean poset $B(3)$ with vertices decorated with the horizontal homologies of a simplicial complex with $3$ vertices, and its flattening to the \"uber chain complex.}\label{fig:cubo}
\end{figure}

For a bi-colouring $\e$ on the simplicial complex $X$, we set $\ell(\e) \coloneqq \sum_{j} \e(v_j) $. We can then define the \emph{$j$-th \"uber chain complex} (with coefficients in $R$) as follows:
\begin{equation}
\UC^{j}(X;R) = \bigoplus_{\ell(\e) = j} \Hh(X,\e;R) \ .
\end{equation}
We now restrict to the case $R=\bF$. Then,  by \cite[Proposition~6.2]{uberhomology}, the map
\begin{equation}
    d^j\coloneqq \sum_{\ell(\e) = j} d_{\e,\e'}\colon \UC^{j}(X;R)\to \UC^{j+1}(X;R)
\end{equation}
is a differential, turning $\left(\UC^{*}(X;\bF), d\right)$ into a triply graded complex. A schematic description of the construction for the \"uber chain complex is presented in~Figure~\ref{fig:cubo}.

\begin{defn}
The \emph{\"uberhomology} $\UH^*(X)$ of a finite and connected simplicial complex $X$ is the homology of the complex $\left(\UC^{*}(X;\bF), d\right)$.
\end{defn}

We refer the reader to \cite[Sections~6,7]{uberhomology} for a detailed construction of the {\"uberhomology}, as well as some examples and computations. 

\begin{rem}\label{rmk:bigrading}
The differential $d$ preserves both the dimension of the simplices and their weight. It follows that the \"uberhomology ``inherits'' two gradings from the horizontal homology, making it a triply-graded homology. As a matter of notation, we will sometimes denote these gradings as $\UH_{i,k}^{j}(X;R)$; the $j$ grading is the homological degree of the \"uberhomology, increasing by $1$ under the action of $d$. The other bi-degree $(i,k)$ denotes the pair consisting of dimension of simplices and weight.
\end{rem}

\subsection{Poset homology}\label{posethomology}
We start by reviewing the poset homology of a \emph{finite} poset~$P$, with coefficients in a functor~$\mathcal{F}$. We remark here that this construction is related, but in general not equivalent, to the classical poset homology (see~\emph{e.g.}~\cite{Wachs}, and \emph{cf.}~Remark~\ref{rmk:posethom}) which is defined as the homology of the associated nerve. We refer to \cite{chandler2019posets, primo} for more general expositions on the topic. 

For a poset~$(P,\triangleleft)$, let $\widetilde{\triangleleft}$ denote the associated covering relation -- i.e.~$x\ \widetilde{\triangleleft} \ y$ if and only if $x\triangleleft y$ and there is no $z$ such that $x\triangleleft z\triangleleft y$. 

We say that $P$ is \emph{ranked} if there is a rank function $\ell\colon P\to \bN$ such that $x\ \widetilde{\triangleleft} \ y$ implies $\ell(y)=\ell(x)+1$. We say that  $P$ is \emph{squared} if, for each triple $x,y,z\in P$ such that $z$ covers $y$ and $y$ covers $x$, then
there is a unique $y' \neq y$ such that $z$ covers $y'$ and $y'$ covers $x$. Such elements $x,y,y',z$, together with their covering relations in $P$, will be called a \emph{square}.
In what follows, we assume all posets to be ranked and squared.

\begin{example}
A Boolean poset is ranked and squared; the rank function is given by the distance~$\ell$ of an element from the empty set (cf.~Remark~\ref{rem:boolean col}).
\end{example}

Observe that a poset can always be regarded as a category:
\begin{rem}\label{rem:posetiscat}
A finite poset $(P,\triangleleft)$ can be  seen as a (small) category $\mathbf{P}$; the set of objects of {\bf P} is the set $P$, and the set of morphisms between $x$ and $y$ contains a single element if and only if $x\triangleleft y$ or $x=y$, and is empty otherwise.
\end{rem}
Functors on the category associated to the poset~$P$ preserve commutative squares; in fact, we have the following:
\begin{rem}\label{rem:squares}
Let ${\bf C}$ be a small category, and $(P,\triangleleft)$ be a poset.
For each $x,z \in P$ there is a unique mapping $f_{x,z}:x\to z$ in the category {\bf P}.
Assume there is a square between $x$ and~$z$; the existence of such a square implies that $f_{x,z}$ factors:
\[ f_{x,z} = f_{y,z} \circ f_{x,y} = f_{y',z} \circ f_{x,y'}.\]
Given a covariant functor $\mathcal{F}: {\bf P} \to {\bf C}$, we must have:
\[\mathcal{F}(f_{y,z}) \circ \mathcal{F}(f_{x,y}) =\mathcal{F}(f_{y,z} \circ f_{x,y}) =\mathcal{F}(f_{x,z})= \mathcal{F}(f_{y',z} \circ  f_{x,y'}) = \mathcal{F}(f_{y',z}) \circ \mathcal{F}(f_{x,y'}) \ . \]
In other words, all functors preserve the commutativity of the squares in $P$.
\end{rem}  

Let $\bZ_2$ be the cyclic group on two elements.

\begin{defn}\label{def:sign_ass}
	A \emph{sign assignment} on a poset $(P,\triangleleft)$ is an assignment of elements $s_{x,y}\in \bZ_2$ to each pair of elements $x,y\in P$ with $x\ \widetilde{\triangleleft}\ y$, such that the equation
	\begin{equation}\label{eq:signassign}
		s_{x,y} + s_{y,z} \equiv s_{x,y'} + s_{y',z} + 1 \mod 2
	\end{equation}
	holds for each square  $x\ \widetilde{\triangleleft}\ y,~y'\ \widetilde{\triangleleft}\ z$. 
\end{defn}

In general, the existence of a sign assignment on a poset~$P$ depends on the topology of a certain topological space associated to~$P$ -- see, \emph{e.g.}~\cite[Section~3.2]{primo}, \cite[Section~5]{Putyra}, or \cite{chandler2019posets}. 
In cases of interest to us, there is always a sign assignment, and the choice of such a sign is thus immaterial.

\begin{rem}\label{rem:Boolean sings}
Any Boolean poset admits a sign assignment, which is unique up to (a suitable notion of) isomorphism -- \emph{cf.}~\cite[Example~3.15]{primo}.
\end{rem}

We can now recall the definition of  poset homology of a poset~$P$ with coefficients in a functor~$\cF$.

Let ${\bf A}$ be an Abelian category -- \emph{e.g.}~the category of left modules on a commutative ring~$R$ --  $P$ a ranked squared poset with rank function~$\ell$, and $s$ a sign assignment on $P$.
Given a covariant functor 
$ \mathcal{F}\colon{\bf P} \to {\bf A}$,
we can define the cochain groups 
\[ C^n_{\mathcal{F}}(P) \coloneqq \bigoplus_{\tiny\begin{matrix}
		x\in P\\
		\ell(x) = n
\end{matrix}}  \mathcal{F}(x),
\]
and the differentials
\[d^{n}=d^{n}_{\mathcal{F}} \coloneqq \sum_{\tiny\begin{matrix}
		x\in P\\
		\ell(x) = n
\end{matrix}} \sum_{\tiny\begin{matrix}
		x'\in P\\
		x \ \widetilde{\triangleleft} \ x'
\end{matrix}} (-1)^{ s(x,x')} \mathcal{F}(x \ \widetilde{\triangleleft} \ x') \ . \]

Note that the differentials $d^n$, and therefore the cochain complexes,  depend \emph{a priori} upon the choice of the sign assignment~$s$. However, in the cases of interest to us -- \emph{i.e.}~for Boolean posets -- this choice does not affect the isomorphism type of the cochain complexes.

\begin{thm}\label{teo: general cohom}
Let ${\bf A}$ be an Abelian category, $P$  be a ranked squared poset, and $s$ be a sign assignment on $P$. Then, for any $n\in\mathbb{N}$ and any functor $\mathcal{F}\colon{\bf P} \to {\bf A}$ we have $d^{n+1} \circ d^{n} \equiv 0$. In particular, $(C^*_{\mathcal{F}}(P), d^*)$ is a cochain complex.
\end{thm}

For a proof of this result, we refer to \cite{chandler2019posets} and \cite[Theorem~3.7]{primo}.

\begin{rem}\label{rmk:posethom}
Poset homology, as described in this section, is related to the classical homology of posets -- defined as the homology of the associated nerve~\cite[Section~1.5]{Wachs}. Indeed, when the poset~$P$ is the face poset of a CW-complex~$X$, and the functor~$\mathcal{F}$ is the constant functor, then the poset homology of $P$ with coefficients in $\mathcal{F}$ agrees with the reduced homology of~$X$  (shifted by~$1$) -- see \cite[Section~6]{secondo}.
Furthermore, when the poset~$P$ is a Boolean poset, the relationship is stronger: the poset homology with coefficients in a functor~$\mathcal{F}$ agrees with the homology (of the associated category) with coefficients in~$\mathcal{F}$ -- now defined as the derived functors of $\mathrm{colim}(\mathcal{F})$~\cite{Gabriel1967CalculusOF} -- by~\cite[Theorem~24]{turner-everitt}. 
\end{rem}

\subsection{\"Uberhomology as a poset homology}\label{sec:uberasposet}

Let $\mathbf{Mod}_R$ be the category of (left)  $R$-modules, over a fixed commutative ring~$R$ with unit.
Note that the category $\mathbf{Mod}_R$ is an Abelian category; in particular, biproducts are given by direct sums of modules.

Let $X$ be a simplicial complex with $|V(X)|=m$. 
By Remark~\ref{rem:boolean col}, we have an identification of each $b\in B(m)$ with a bi-colouring $\e_{b}$ on $X$. Consider the category $\mathbf{B}(m)$ associated to $B(m)$, as in  Remark~\ref{rem:posetiscat}. Then, we can regard the decoration provided in Section~\ref{sec:uberhom} as a functor 
\[
\cH\colon \mathbf{B}(m)\to \mathbf{Mod}_R
\]
defined as  $\cH(b) \coloneqq \Hh(X,\e_b;R)$ on objects, and as $\cH(b \  \widetilde{\triangleleft} \  b')\coloneqq  d_{\e_b,\e_{b'}}$ for each covering relation $b \  \widetilde{\triangleleft} \  b'$ of $B(m)$. The extension to other morphisms of  $\mathbf{B}(m)$ is obtained by compositions. 
Furthermore, the assignment so described does indeed define a functor by \cite[Section~3]{chandler2019posets}.

\begin{rem}
As the horizontal homology is bi-graded, the functor $\cH$ lands in the subcategory of $\mathbf{Mod}_R$ given by bi-graded modules.
\end{rem}

We are now ready to identify the \"uberhomology with a poset homology on $B(m)$. 

\begin{prop}
Let $X$ be a finite connected simplicial complex with $|V(m)|=m$.
Then, the poset homology of $B(m)$ with coefficients in the functor $\mathcal{H}\colon \mathbf{B}(m)\to \mathbf{Mod}_{\bF}$ agrees with the \"uberhomology of  $X$ with coefficients in $\bF$.
\end{prop}

\begin{proof}
It is enough to write down the definition of $(C^*_{\mathcal{H}}(B(m)), d^*)$, and compare it with the definition of \"uberhomology. The $n$-th cochain group is given by 
\[ C^i_{\mathcal{H}}(B(m)) \coloneqq \bigoplus_{\tiny\begin{matrix}
		b\in B(m)\\
		\ell(\e_b) = i
\end{matrix}}  \mathcal{H}(b) = 
\bigoplus_{\tiny\begin{matrix}
		b\in B(m)\\
		\ell(\e_b) = i
\end{matrix}}  \Hh(X,\e_b;\bF) = \UC^{i}(X;\bF).
\]
Similarly, the differential is given by
\[d^{i}_{\cH} \coloneqq \sum_{\tiny\begin{matrix}
		b\in B(m)\\
		\ell(\e_b) = i
\end{matrix}} \quad  \sum_{\tiny\begin{matrix}
		b'\in B(m)\\
		b \ \widetilde{\triangleleft} \ b'
\end{matrix}} \mathcal{H}(b\ \widetilde{\triangleleft} \ b') \ \overset{(\ast)}{=}
\sum_{\ell(\e) = i}\quad  \sum_{\ell(\e') = i + 1}  d_{\e,\e'} = d_i, \]
where in $(\ast)$ we used that by definition $d_{\e_b,\e_{b'}} = 0$ if $b$ is not covered by $b'$.
\end{proof}

This alternative description of the \"uberhomology as a poset homology allows us, by Theorem~\ref{teo: general cohom}, to extend its definition to any ring of coefficients~$R$. In fact, once a sign assignment~$s$ on $B(m)$ is chosen, the \"uberhomology with coefficients in $R$ is defined as the poset homology of~$B(m)$ with coefficients in the functor
$\mathcal{H}\colon \mathbf{B}(m)\to \mathbf{Mod}_{R}$. 
Note that, up to isomorphism, the result does not depend on the chosen sign assignment by Remark~\ref{rem:Boolean sings}. 

The following corollary describes the functoriality of \"uberhomology with respect to such coefficients:

\begin{cor}
For each finite simplicial complex~$X$, its \"uberhomology defines a functor
\[
\UH^*(X;-)\colon \mathbf{Ring}\to \mathbf{grAb}.
\]
That is, for each ring homomorphism $\phi\colon R\to S$ there is an induced map $\phi^*\colon \UH^*(X;R)\to \UH^*(X;S)$ of graded Abelian groups.
\end{cor}

\begin{proof}
A homomorphism of rings $\phi\colon R\to S$  induces a natural transformation between the functors $\mathcal{H}_R\colon \mathbf{B}(m)\to \mathbf{Mod}_{R}$ and $\mathcal{H}_S\colon \mathbf{B}(m)\to \mathbf{Mod}_{S}$ by extension of scalars. The result is then a consequence of \cite[Corollary~7.15]{chandler2019posets}.
\end{proof}

\section{Functoriality of \"uberhomology}\label{sec:functoriality}

In this section we investigate \"uberhomology's functoriality with respect to certain simplicial maps.
Recall that a simplicial map is a map between simplicial complexes such that the image of the vertices of any simplex spans a simplex. 

\begin{defn}\label{def:col map}
A \emph{coloured map} $\psi\colon(X,\e_X) \to (Y,\e_{Y}) $ is a simplicial map $\psi:X\to Y$ which, for each $y \in \psi (V(X)) \subseteq V(Y)$, satisfies the following two conditions: 
\begin{enumerate}[label={\rm (\arabic*)}]
    \item there is at most one $x\in V(X)$ such that {$\e_{X}(x) = 1$}  and $\psi(x)=y$;
    \item {$\e_Y (y) = 0$}  if and only if there is no $x\in V(X)$ such that {$\e_{X}(x) = 1$}  and $\psi(x)=y$.
\end{enumerate}
\end{defn}

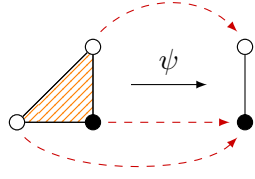
\begin{figure}[ht]
    \centering
    \begin{tikzpicture}

    \draw[pattern = north east lines, pattern color = orange] (0,0) --  (1,1)-- (1,0)  -- (0,0);
    \draw[fill] (0,0) circle (.1) --  (1,1) circle (.1) -- (1,0) circle (.1) -- (0,0);
    \draw[fill, white] (1,1) circle (.1);
    \draw[] (1,1) circle (.1);
    
    \draw[fill, white] (0,0) circle (.1);
    \draw[] (0,0) circle (.1);
    
    \begin{scope}[shift={+(2,0)}]
    \draw[fill] (1,0) circle (.1) --  (1,1) circle (.1);
    \draw[fill, white] (1,1) circle (.1);
    \draw[] (1,1) circle (.1);
    \end{scope}
    
    \draw[-latex] (1.5,.5) -- (2.5,0.5);
    \node[above] at (2,.5) {$\psi$};
    
    \draw[bunired, dashed, -latex] (1.1,1.2) .. controls +(.5,.5) and +(-.5,.5) .. (2.9,1.2);
    \draw[bunired, dashed, -latex] (1.2,0) -- (2.8,0);
    \draw[bunired, dashed, -latex] (0.1,-.2) .. controls +(.5,-.5) and +(-.5,-.5) .. (2.9,-.2);
    \end{tikzpicture}
    \caption{A coloured map $\psi$; the collapse of a coloured $2$-simplex on one of its edges.}
    \label{fig:colour-compatible simplicial not inj}
\end{figure}

Figure~\ref{fig:colour-compatible simplicial not inj} gives an example of a coloured map of simplicial complexes. The identity is always a coloured map between $(X,\e_X)$ and itself. Coloured maps are closed under composition:

\begin{lem}
The composition of two composable coloured maps is a coloured map.
\end{lem}

\begin{proof}
Composition of simplicial maps is simplicial. It is also straightforward to check that the conditions defining a coloured map are satisfied.
\end{proof}

\begin{rem}\label{rem:inj col maps}
Observe that, given a coloured map $\psi\colon(X,\e_X) \to (Y,\e_{Y})$, then {$\e_{Y}(v) = 0$} for each $v\in V(Y)\setminus \psi(V(X))$.
\end{rem}

Coloured simplicial maps are compatible with the  weight of  bi-coloured simplicial complexes, as defined in Equation~\eqref{eq:weight}:

\begin{lem}\label{lem:w pres}
Let $\sigma$ be a simplex of $X$,  and $\psi\colon(X,\e_X) \to (Y,\e_{Y})$  be a coloured map. Then,
\[
w_{\e_X}(\sigma ) \geq w_{\e_Y} (\psi(\sigma) ) \ ,
\]
and equality holds if $\psi$ is injective.
That is, injective coloured simplicial maps preserve weights. 
\end{lem}

\begin{proof}
Note that a coloured map $\psi\colon(X,\e_X) \to (Y,\e_{Y})$ preserves the sum of the colourings on the vertices of each simplex, that is the equality
\[ \sum_{v\in V(\sigma)} \e_X(v)  = \sum_{x\in V(\psi(\sigma))} \e_Y (x) \]
holds for each $\sigma \subseteq X$. The image $\psi(\sigma)$ is a simplex in $Y$, of possibly lower dimension, with the same number of $1$-coloured vertices as~$\sigma$.
The statement now follows from the definition of weight.
\end{proof}

We observe here that a generic coloured simplicial map may not induce a chain map between the associated $\e$-horizontal chain complexes:

\begin{rem}\label{rem:counterx}
Fix a colouring $\e_{\Delta^1}$ on the standard $1$-simplex $\Delta^1$ and a colouring $\e_{\Delta^0}$ on the standard $0$-simplex $\Delta^0$. 
Let $\psi\colon\Delta^1\to \Delta^0$ be the simplicial map obtained by collapsing~$\Delta^1$ to~$\Delta^0$. Then, $\psi$ uniquely determines a morphism $\overline{\psi} \colon C(\Delta^1,\e_{\Delta^1};R)\to C(\Delta^0,\e_{\Delta^0};R)$ of $R$-modules, since the underlying $R$-module of the horizontal chain complex $C(X,\e;R)$ of a coloured simplicial complex~$(X,\e)$ is isomorphic to the $R$-module of the standard simplicial chain complex associated to $X$. However, $\overline{\psi}$ does not extend to a chain map with respect to the horizontal differentials. 
\end{rem}

Note that,  when restricting to injective coloured simplicial maps, we get graded maps of complexes with respect to the weight grading; more specifically, for any coefficient ring $R$, we have the following:

\begin{prop}\label{prop:inj col}
An injective coloured simplicial map $\psi\colon (X,\varepsilon_X) \to (Y,\varepsilon_{Y}) $ induces a chain map between the associated horizontal chain complexes
\[ \overline{\psi} \colon ({\rm C}(X,\e_X;R),\partial_h ) \longrightarrow ({\rm C}(Y,\e_Y;R),\partial_h) \ , \]
which is graded with respect to the grading induced by the weights $w_{\e_X}$ and $w_{\e_Y}$.
\end{prop}

\begin{proof}
Simplicial maps induce morphisms of chain complexes. Furthermore, for an injective map, the weights are preserved by Lemma~\ref{lem:w pres}, making the induced morphism of chain complexes  a graded one.
\end{proof}

As a consequence of the above proposition, we obtain functoriality of the $\e$-horizontal homology $\Hh(X,\e;R)$ with respect to injective coloured simplicial maps.

We now restrict to injective simplicial maps of simplicial complexes.

\begin{thm}\label{thm:functo simpl ,map}
For each commutative ring~$R$, the \"uberhomology is a functor
\[
\UH(-;R)\colon \mathbf{RegSimpl}\to \mathbf{Mod}_R.
\]
from the category $\mathbf{RegSimpl}$ of simplicial complexes and injective simplicial maps.
\end{thm}

\begin{proof}
Let $\varphi\colon X\to Y$ be an injective simplicial map, and assume $|V(X)|=m$ and $|V(Y)|=n$. The Boolean poset $B(m)$ can be seen as a sub-poset of $B(n)$; consider the functors $
\cH_X\colon \mathbf{B}(m)\to \mathbf{Mod}_R$ and $
\cH_Y\colon \mathbf{B}(n)\to \mathbf{Mod}_R
$
as described in Section~\ref{sec:uberasposet}. Note that the functor $\cH_X$ can be uniquely extended to a functor~$\widetilde{\cH}_X$ on the category $\mathbf{B}(n)$ (via the zero extension); this extension provides a natural transformation $\eta\colon \widetilde{\cH}_X \Rightarrow  \cH_Y$. The result now follows by \cite[Corollary~7.15]{chandler2019posets}.
\end{proof}

The functoriality of Theorem~\ref{thm:functo simpl ,map} cannot be extended to the whole category $\mathbf{Simpl}$ of simplicial complexes with the same methods -- see also Remark~\ref{rem:counterx}. With a different approach, one may switch the  roles of $0$ and $1$ in Definition~\ref{def:col map}, obtaining the analogue of Proposition~\ref{prop:inj col} for \emph{all} coloured (with respect to this new definition) simplicial maps. However, in such case, the functoriality -- again only for injective simplicial maps -- of Theorem~\ref{thm:functo simpl ,map}
is only true up to a shift in the degree.

\section{Bold homology and Dominating sets}\label{sec: bold homology and dom}

In this section we review the definition of bold homology, and investigate some of its properties.

\subsection{Bold homology}

The \"uberhomology introduced in Section~\ref{sec:uberhom} is a homology theory for simplicial complexes; when restricting to simple graphs, \emph{i.e.}~$1$-dimensional simplicial complexes, we can consider what in  \cite[Section~8]{uberhomology} has been denoted by $\bH^0$. For notational convenience in what follows we will drop the index~$0$ from the notation. This is a singly-graded homology theory,  consisting of the bidegree $(0,0)$ part of the \"uberhomology, and called here \emph{bold homology}. We start by recalling  its alternative description in terms of connected subgraphs.

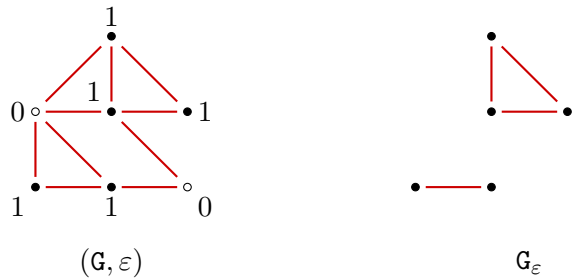
\begin{figure}[ht]
\centering
\begin{tikzpicture}
    \node (a) at (0,0) {};
    \node (b) at (1,0) {};
    \node (c) at (0,1) {};
    \node (d) at (1,-1) {};
    \node (e) at (-1,0) {};
    \node (f) at (0,-1) {};
    \node (g) at (-1,-1) {};

    \node[above left]   at (0,0) {1};
    \node[right]   at (1,0) {1};
    \node[above]   at (0,1) {1};
    \node[below right]   at (1,-1) {0};
    \node[left]   at (-1,0) {0};
    \node[below]   at (0,-1) {1};
    \node[below left]   at (-1,-1) {1};
    
    \draw[fill] (a) circle (.05);
    \draw[fill] (b) circle (.05);
    \draw[fill] (c) circle (.05);
    \draw[    ] (d) circle (.05);
    \draw[    ] (e) circle (.05);
    \draw[fill] (f) circle (.05);
    \draw[fill] (g) circle (.05);

    \draw[thick , bunired] (a) -- (b);
    \draw[thick , bunired] (a) -- (c);
    \draw[thick , bunired] (a) -- (d);
    \draw[thick , bunired] (a) -- (e);
    \draw[thick , bunired] (c) -- (b);
    \draw[thick , bunired] (c) -- (e);
    \draw[thick , bunired] (f) -- (e);
    \draw[thick , bunired] (f) -- (d);
    \draw[thick , bunired] (g) -- (e);
    \draw[thick , bunired] (f) -- (g);
    
    \node at (0,-2) {$(\tG,\e)$};
    
    \begin{scope}[shift = {+(5,0)}]
    \node (a) at (0,0) {};
    \node (b) at (1,0) {};
    \node (c) at (0,1) {};
    \node (f) at (0,-1) {};
    \node (g) at (-1,-1) {};

    \draw[fill] (a) circle (.05);
    \draw[fill] (b) circle (.05);
    \draw[fill] (c) circle (.05);
    \draw[fill] (f) circle (.05);

    \draw[fill] (g) circle (.05);
    
    \draw[thick , bunired] (a) -- (b);
    \draw[thick , bunired] (a) -- (c);
    \draw[thick , bunired] (f) -- (g);
    \draw[thick , bunired] (c) -- (b);
    
    \node at (0.5,-2) {$\tG_\e$};
    \end{scope}
\end{tikzpicture}
\caption{A bi-coloured graph $(\tG,\e)$, and the subgraph $\tG_\e$ determined by its colouring~$\e$.}\label{fig:col subgraph}
\end{figure}

Let $\tG$ be a connected simple graph.
Denote by $\cl(\tG)$ the set of bi-colourings of $\tG$, and set
\[\cl(\tG;i) \coloneqq \Big\{   \e \in \cl(\tG) : \sum_{v \in V(\tG)} \e (v) = i \Big\}\ .\]
\begin{rem}\label{rem:col and Bool}
The set $\cl(\tG)$ has a natural poset structure induced by the $\e$-colourings and a fixed order on $V(\tG)$. We can identify this poset with a Boolean $B(|V(\tG)|)$ -- \emph{cf.}~Remark~\ref{rem:boolean col}. The choice of the ordering is immaterial, up to poset isomorphism. From now on, when identifying the set~$\cl(\tG)$ with a Boolean poset, we always implicitly assume a choice of an ordering on~$V(\tG)$.
\end{rem}
Each bi-colouring $\e \in \cl(\tG)$ determines a possibly disconnected subgraph $\tG_{\e} \subseteq \tG$ whose vertices are given by $V(\tG_\e) =\{ v \in V(\tG) \mid \e(v) = 1 \}$, and containing all the edges in $\tG$ connecting vertices in $V(\tG_\e)$ -- see Figure~\ref{fig:col subgraph} for an example.

The \emph{$i$-th bold chain group} (with coefficients in $R$) is
\[ {\rm C}\bH_i({\tt G};R) = \bigoplus_{\e \in \cl(\tG;i)} \bigoplus_{x \in \pi_0(\tG_\e)} R \langle x \rangle, \]
where $\pi_0$ denotes the set of connected components.
Note that for $i=0$, $\cl(\tG,0)$ contains only the trivial (all $0$) colouring $\e_0$. The graph associated to $\tG_{\e_0}$ is the empty graph, therefore we can set ${\rm C}\bH_0({\tt G};R) = (0)$.

When clear from the context, we identify a connected component $x$ in $\pi_{0}(\tG_\e)$ with the corresponding subgraph of $\tG$. Given $x\in \pi_{0}(\tG_\e)$ and $x'\in \pi_{0}(\tG_{\e'})$, we write~$x\prec x'$ if and only if $ \e \ \widetilde{\triangleleft} \ \e'$ (under the identification of $\cl(\tG)$ with a suitable Boolean poset -- \emph{cf.}~Remark~\ref{rem:col and Bool}) and~$V(x) \subseteq V(y)$.

Now, define the map of $R$-modules
\[
\begin{matrix}
d^i \colon  & {\rm C}\bH_i(\tG)& \longrightarrow &  {\rm C}\bH_{i+1}(\tG) \\  
& x &  \longmapsto  &\sum_{x\prec  y} s(\e,\e')y \\
\end{matrix}
\]
on generators $x\in \pi_0(\tG_\e)$, for $\e \in \cl(\tG,i)$, and extended by $R$-linearity.
It turns out that~$d$ defines a differential on ${\rm C}\bH_*(\tG;R)$. This is implicit from the fact that $d$  is the component of bi-degree $(0,0)$ of the \"uberhomology differential. We provide a direct proof of this fact.

\begin{lem}
Let $\tG$ be a simple graph. Then, $({\rm C}\bH_*(\tG;R), d^*)$ is a chain complex.
\end{lem}

\begin{proof}
In order to prove that $d$ squares to $0$, consider the composition
\[ d^{i+1} \circ d^{i} (x) = \sum_{x\prec  y} \ \sum_{y\prec  z} s(\e,\e') s(\e',\e'') \ z \ . \]
Note that each $z\in \pi_0(\tG_{\e''})$ appears exactly twice since the poset of colourings is squared. More precisely, there are exactly two colourings $\e'_1$ and $\e'_2$ such that $\e\ \widetilde{\triangleleft} \ \e'_1$, and $\e'_2  \widetilde{\triangleleft}\ \e'' $. By definition of sign assignments, we have $$s(\e,\e'_1)s(\e'_1,\e'') = - s(\e,\e')s(\e',\e'')$$ and the statement follows.
\end{proof}

The homology of $({\rm C}\bH_*(\tG;R), d^*)$ is denoted by $\mathbb{H}_*(\tG; R)$, and referred to as \emph{bold homology}. We simply write $\mathbb{H}_*(\tG)$ when the base ring $R$ is clear from the context.

We can now analyse the bold homology groups in their lowest degrees; first note that by definition ${\rm C}\bH_0(\tG) = (0)$. In particular, $\bH_{0}(\tG)$ is always trivial.
We also have a complete characterisation of the first bold homology group $\mathbb{H}_1$:

\begin{prop}\label{prop: complete graphs}
Let $\tG$ be a simple {connected} graph. Then, 
$\bH_1(\tG) \neq 0$ if and only if $\tG \cong \tK_m$ for  some $m\ge0$.
\end{prop}
We sketch a proof here, and give a full proof in Section~\ref{sec: retrction} using the techniques developed therein.
\begin{proof}[Sketch of proof]
If $\tG$ is a complete graph, then $\bH_1(\tG) \cong \bF$ by \cite[Proposition~8.1]{uberhomology}; note that this also holds on any field $\bK$. 
Complete graphs are the only graphs such that the subgraphs induced from a $1$-colouring on precisely two vertices are always connected. Therefore, if $\tG$ is not a complete graph, there exists at least a pair of vertices $v_1,v_2 \in V(\tG)$ such that the induced graph is disconnected. We use this to show that the differential~$d^1$ is injective.

First note that if a linear combination of generators $\alpha_1 v_{i_1} + \ldots +\alpha_k v_{i_k}$ contains $v_1$ or $v_2$, then its image $d^1(\alpha_1 v_{i_1} + \ldots +\alpha_k v_{i_k})$ can not be trivial (see Figure~\ref{fig:onlycomplete}). Note also that every connected generator in degree two is in the image of exactly two connected components (\emph{i.e.}~two $1$-coloured vertices). If a vertex $v$ is not connected to all the other vertices of $\tG$, we can always find a $v'$ such that $v \cup v'$ induces a disconnected subgraph. It follows that, if $\alpha_1 v_{i_1} + \ldots +\alpha_k v_{i_k}$ is a cycle, then each $v_{i_j}$ must be connected with every other vertex in $\tG$.

Now, given $v\in V(\tG)\setminus \{ v_1 \}$, connected with all the other vertices in $\tG$, we have that $\{ v_1, v\} $ appears in $d^1 v$ with coefficient $\pm 1$. Hence, to cancel out this contribution, $v_1$ must appear in each cycle featuring $v$, which is a contradiction. 
\begin{figure}[ht]

\labellist
\pinlabel $4$ at 48 230
\pinlabel $1$ at 48 360
\pinlabel $2$ at -16 295
\pinlabel $3$ at 109 295
\endlabellist
\centering
\includegraphics[width = 6cm]{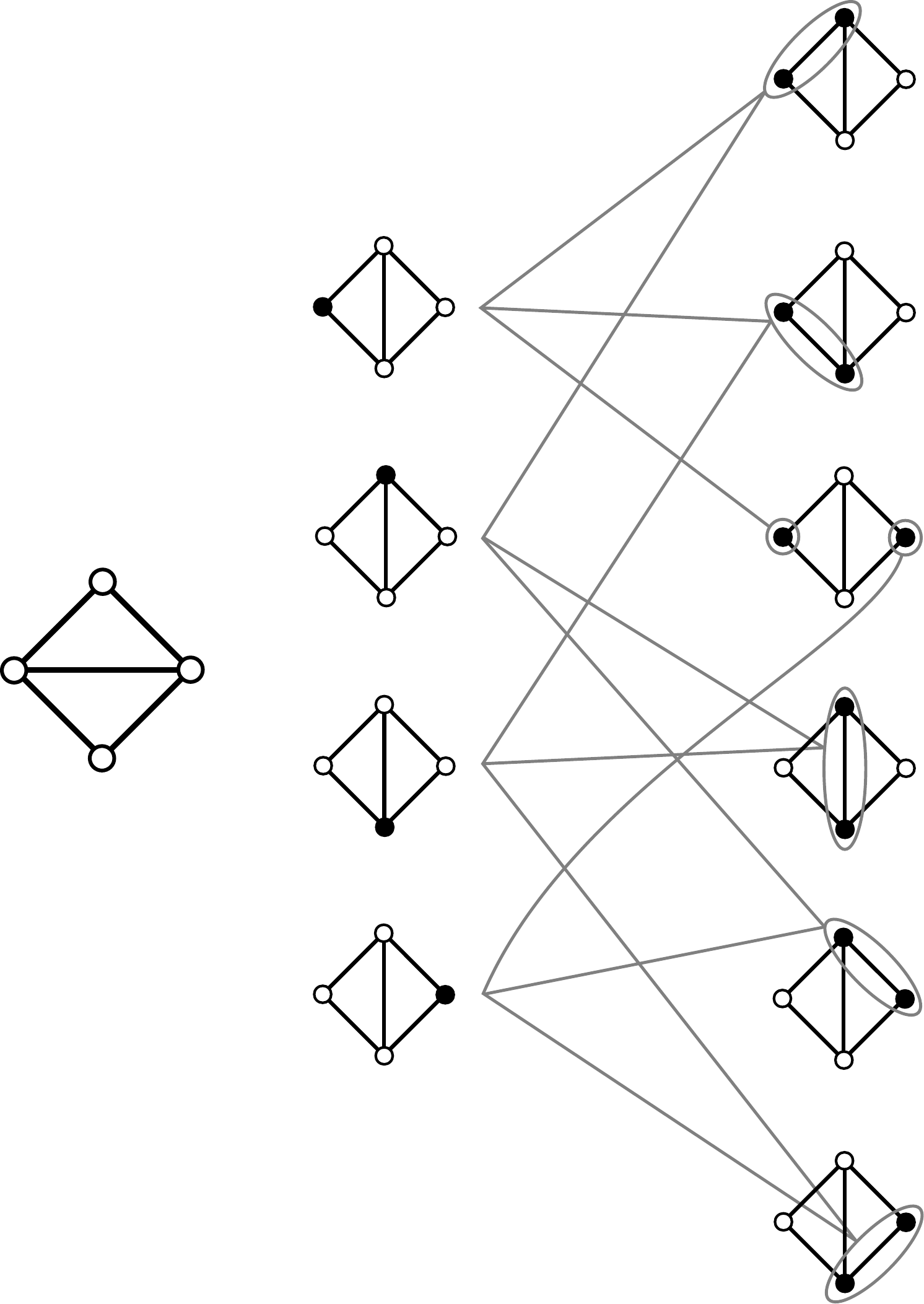}
\caption{A graph $\tG$ on $4$ vertices, together with its bold chain complex in the first two degrees. Since $\tG$ is not a complete graph, we can see two connected components appearing in degree $2$. In the case at hand, these correspond to the $\e$-colouring $(1,0,0,1)$.}
\label{fig:onlycomplete}
\end{figure}
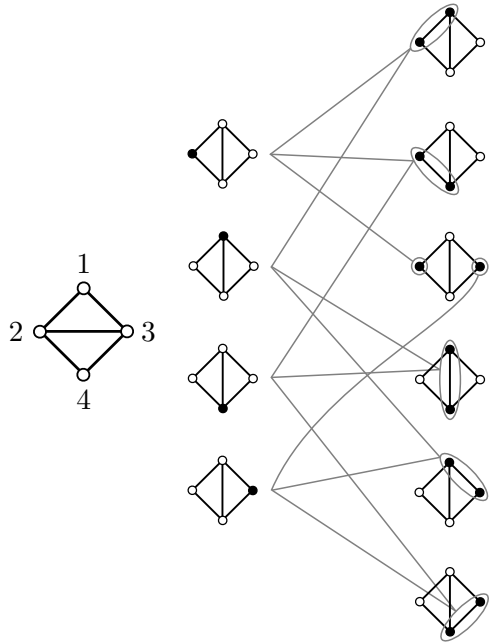
\end{proof}

\subsection{Dominating sets and Domination polynomials}\label{sec: dom sets and poly}

In this section we introduce some basic notions related to graphs and dominating sets that will be used throughout the rest of the paper, and show a categorification-like relation with bold homology.

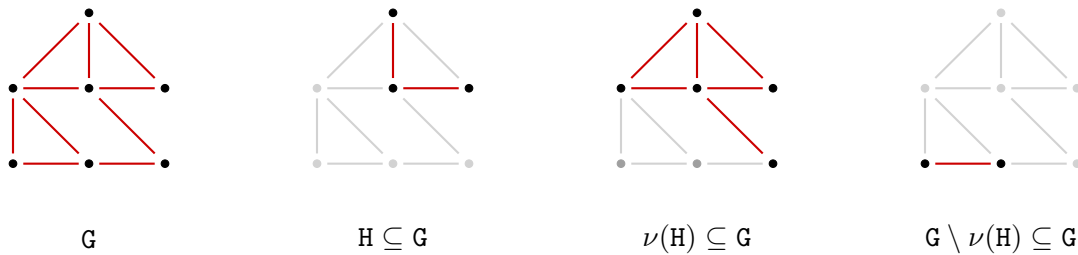
\begin{figure}[ht]
\centering
\begin{tikzpicture}
    \node (a) at (0,0) {};
    \node (b) at (1,0) {};
    \node (c) at (0,1) {};
    \node (d) at (1,-1) {};
    \node (e) at (-1,0) {};
    \node (f) at (0,-1) {};
    \node (g) at (-1,-1) {};

    \draw[fill] (a) circle (.05);
    \draw[fill] (b) circle (.05);
    \draw[fill] (c) circle (.05);
    \draw[fill] (d) circle (.05);
    \draw[fill] (e) circle (.05);
    \draw[fill] (f) circle (.05);
    \draw[fill] (g) circle (.05);

    \draw[thick , bunired] (a) -- (b);
    \draw[thick , bunired] (a) -- (c);
    \draw[thick , bunired] (a) -- (d);
    \draw[thick , bunired] (a) -- (e);
    \draw[thick , bunired] (c) -- (b);
    \draw[thick , bunired] (c) -- (e);
    \draw[thick , bunired] (f) -- (e);
    \draw[thick , bunired] (f) -- (d);
    \draw[thick , bunired] (g) -- (e);
    \draw[thick , bunired] (f) -- (g);
    
    \node at (0,-2) {$\tG$};
    
    \begin{scope}[shift = {+(4,0)}]
    \node (a) at (0,0) {};
    \node (b) at (1,0) {};
    \node (c) at (0,1) {};
    \node (d) at (1,-1) {};
    \node (e) at (-1,0) {};
    \node (f) at (0,-1) {};
    \node (g) at (-1,-1) {};

    \draw[fill] (a) circle (.05);
    \draw[fill] (b) circle (.05);
    \draw[fill] (c) circle (.05);
    \draw[fill, gray,  opacity =.35] (d) circle (.05);
    \draw[fill, gray,  opacity =.35] (e) circle (.05);
    \draw[fill, gray,  opacity =.35] (f) circle (.05);
    \draw[fill, gray,  opacity =.35] (g) circle (.05);

    \draw[thick , bunired] (a) -- (b);
    \draw[thick , bunired] (a) -- (c);
    \draw[thick , gray,  opacity =.35] (a) -- (d);
    \draw[thick , gray,  opacity =.35] (a) -- (e);
    \draw[thick , gray,  opacity =.35] (c) -- (b);
    \draw[thick , gray,  opacity =.35] (c) -- (e);
    \draw[thick , gray,  opacity =.35] (f) -- (e);
    \draw[thick , gray,  opacity =.35] (f) -- (d);
    \draw[thick , gray,  opacity =.35] (g) -- (e);
    \draw[thick , gray,  opacity =.35] (f) -- (g);
    
    \node at (0,-2) {$\tH\subseteq \tG$};
    \end{scope}
    
\begin{scope}[shift = {+(8,0)}]
    \node (a) at (0,0) {};
    \node (b) at (1,0) {};
    \node (c) at (0,1) {};
    \node (d) at (1,-1) {};
    \node (e) at (-1,0) {};
    \node (f) at (0,-1) {};
    \node (g) at (-1,-1) {};

    \draw[fill] (a) circle (.05);
    \draw[fill] (b) circle (.05);
    \draw[fill] (c) circle (.05);
    \draw[fill] (d) circle (.05);
    \draw[fill] (e) circle (.05);
    \draw[fill, gray,  opacity =.75] (f) circle (.05);
    \draw[fill, gray,  opacity =.75] (g) circle (.05);

    \draw[thick , bunired] (a) -- (b);
    \draw[thick , bunired] (a) -- (c);
    \draw[thick , bunired] (a) -- (d);
    \draw[thick , bunired] (a) -- (e);
    \draw[thick , bunired] (c) -- (b);
    \draw[thick , bunired] (c) -- (e);
    \draw[thick , gray,  opacity =.35] (f) -- (e);
    \draw[thick , gray,  opacity =.35] (f) -- (d);
    \draw[thick , gray,  opacity =.35] (g) -- (e);
    \draw[thick , gray,  opacity =.35] (f) -- (g);
    
    \node at (0,-2) {$\nu(\tH)\subseteq \tG$};
    \end{scope}
        
\begin{scope}[shift = {+(12,0)}]
    \node (a) at (0,0) {};
    \node (b) at (1,0) {};
    \node (c) at (0,1) {};
    \node (d) at (1,-1) {};
    \node (e) at (-1,0) {};
    \node (f) at (0,-1) {};
    \node (g) at (-1,-1) {};

    \draw[fill, gray,  opacity =.35] (a) circle (.05);
    \draw[fill, gray,  opacity =.35] (b) circle (.05);
    \draw[fill, gray,  opacity =.35] (c) circle (.05);
    \draw[fill, gray,  opacity =.35] (d) circle (.05);
    \draw[fill, gray,  opacity =.35] (e) circle (.05);
    \draw[fill] (f) circle (.05);
    \draw[fill] (g) circle (.05);

    \draw[thick, gray,  opacity =.35] (a) -- (b);
    \draw[thick, gray,  opacity =.35] (a) -- (c);
    \draw[thick, gray,  opacity =.35] (a) -- (d);
    \draw[thick, gray,  opacity =.35] (a) -- (e);
    \draw[thick, gray,  opacity =.35] (c) -- (b);
    \draw[thick, gray,  opacity =.35] (c) -- (e);
    \draw[thick , gray,  opacity =.35] (f) -- (e);
    \draw[thick , gray,  opacity =.35] (f) -- (d);
    \draw[thick , gray,  opacity =.35] (g) -- (e);
    \draw[thick , bunired] (f) -- (g);
    
    \node at (0,-2) {$\tG\setminus \nu(\tH)\subseteq \tG$};
    \end{scope}
\end{tikzpicture}
\caption{A graph $\tG$, a (non-dominating) connected subgraph $\tH\subseteq \tG$, the neighbourhood $\nu(\tH)$ of $\tH$ in $\tG$,  and its complement $\tG\setminus \nu(\tH)$.}\label{fig:col subgraph and ngbh}
\end{figure}

Let $\tG$ be a simple graph. 
For a subgraph $\tH\subseteq \tG$ we denote by $\nu(\tH)$ the $1$-neighbourhood  of $\tH$ in $\tG$, \emph{i.e.}~the subgraph induced by $V(\tH)$ and by all the vertices sharing an edge with an element of $V(\tH)$. 
The graph $\tG\setminus \nu(\tH)$ is the subgraph of $\tG$ induced by the vertices $V(\tG)\setminus V(\nu(\tH))$ (see  Figure~\ref{fig:col subgraph and ngbh}).
With the poset structure induced by $\cl(\tG)$, the set 
\[\cl(\tG;\tH) = \{ \e \in \cl(\tG)\mid \e(v) =1 \text{ if }v\in V(\tH)\text{, and   }\e(v) = 0 \text{ if }v\in V(\nu(\tH) \setminus \tH)\}\]
is a Boolean sub-poset of $\cl(\tG)$ isomorphic to $B(| V(\tG \setminus \nu(\tH))|)$ -- \emph{cf.}~Remark~\ref{rem:col and Bool}. 

Dominating set have been extensively studied \cite{allan1978domination, alikhani2009introduction,haynes2013fundamentals}; this is mainly due to the NP-complete status of the problem of finding all dominating sets of a graph~\cite{crescenzi1995compendium}, as well as their relationship with the study of networks (see \emph{e.g.}~\cite{wu1999calculating}). We recall here the definition:

\begin{defn}
Given a simple graph $\tG$, a subset $D \subseteq V(\tG)$ is said to be \emph{dominating} if every vertex in $V(\tG) \setminus D$ is adjacent to some member of $D$.
\end{defn}

Equivalently, a subgraph $\tH \subseteq\tG$ is dominating if and only if its $1$-neighbourhood~$\nu(\tH)$ is $\tG$.
The \emph{dominating polynomial} $D_{\tG}(x) \in \mathbb{Z}[x]$ is the polynomial  
\[ D_{\tG}(x) = \sum_{k = \gamma(\tG)}^{|V(\tG)|} d_{\tG}(k) x^k\ ,\]
where $d_{\tG}(k)$ is the number of dominating sets in $\tG$ composed by exactly $k$ vertices, and $\gamma(\tG)$ is the \emph{domination number}, \emph{i.e.} the minimal size of a dominating set in $\tG$.

We are interested in the related notion of  connected dominating sets (see \emph{e.g.}~\cite{sampathkumar1979connected}, or \cite{du2012connected}, for an overview of the applications).

\begin{defn}\label{def:conn dom}
A dominating set is called \emph{connected} if its induced graph is connected. 
The \emph{connected domination polynomial} is defined as follows
\[ D^c_{\tG}(x) = \sum_{k = \gamma^c(\tG)}^{|V(\tG)|} d^c_{\tG}(k) x^k,
\]
where $d^c_{\tG}(k)$ is the number of connected dominating sets with exactly $k$ vertices, and $\gamma^c(\tG)$ the \emph{connected domination number} -- that is the minimal size of a connected dominating set in $\tG$.
\end{defn}

The proof of the following result, relating bold homology and connected dominating sets, was suggested  by F.~Petrov~\cite{russo}.
\begin{prop}\label{prop:russo}
The Euler characteristic $\chi(\bH(\tG))$ of $\bH(\tG)$ coincides with $D^c_{\tG}(-1)$.
\end{prop}
\begin{proof}
Let $\tG$ be a simple graph with $n$ vertices.
Then by definition, we have:  
\begin{equation}\label{eqn:russo}
\chi(\bH(\tG)) = \sum_{k = 1}^n (-1)^k \sum_{\e \in \cl(\tG;k)} \vert \pi_0(G_\e) \vert   = \sum_{k = 1}^n (-1)^k \sum_{\e \in \cl(\tG;k)} \sum_{\tH \in \pi_0(\tG_\e)} 1 \ .
\end{equation}
We can rearrange the summands in the right-hand side of \eqref{eqn:russo}; that is, we count (with sign) how many times a connected subgraph $\tH \subseteq \tG$, induced by its vertices, appears as a component of some $\tG_\e$. After this re-arrangement, we obtain the equality
\[
\chi(\bH(\tG)) = \sum_{\tH} \sum_{\e\in \cl(\tG;\tH)} (-1)^{\ell(\e)} \ .
\]
where $\tH$ ranges among all connected subgraphs of $\tG$ induced by their vertices.
As pointed out above at the beginning of this section, $\cl(\tG;\tH)$ is  a Boolean sub-set of $\cl(\tG)$, and it is easy to see that the sum~$\sum_{\e\in \cl(\tG;\tH)} (-1)^{\ell(\e)}$ is zero, unless $\cl(\tG;\tH)$ contains a single element. This happens if and only if $\tG\setminus \nu(\tH)$ is empty, hence when $\tH$ is dominating (note that $\tH$ is connected by definition). In which case, we have $(-1)^{\ell(\e)} = (-1)^{|V(\tH)|}$, where $\e$ is the unique colouring in $\cl(\tG;\tH)$, and the statement follows.
\end{proof}

Connected domination polynomials have been computed in some cases. This allows us to easily compute the Euler characteristic of $\bH$, as shown in the following example.
\begin{example}
Let  $\tP$ be the Petersen graph.
Its connected domination polynomial was computed in \cite{mojdeh2018connected} to be 
\[D^c_{\tP}(x) = x^{10} + 10x^{9} + 45x^{8} + 110x^7 + 135x^6 + 72x^5 + 10x^4.\]
It follows that $\chi(\bH({\tP})) = D^c_{{\tP}}(-1) = -1$.
This is coherent with our computations; we explicitly computed the bold homology using a computer program, obtaining
\[ \bH_i(\tP;\bF) \cong \begin{cases} \bF & \text{if }i=5, \\ 0 & \text{otherwise}\end{cases}\]
where $\bF$ is the field with two elements.
\end{example}

We are now ready to give a proof of Theorem~\ref{thm:categorification}.

\begin{proof}[Proof of Theorem~\ref{thm:categorification}]
By Theorem~\ref{thm:functo simpl ,map}, the bold homology is a functor from the category of graphs and injective maps to the category  of $R$-modules. This, together with  Proposition~\ref{prop:russo}, provides the needed categorification.
\end{proof}

\subsection{Retraction onto dominating complex}\label{sec: retrction}
In this section we provide a ``categorified version" of Proposition~\ref{prop:russo};
more precisely, we prove the existence of a retraction of the bold chain complex~$ {\rm C}\mathbb{H}({\tt G})$ on the subcomplex~$ {\rm D}\mathbb{H}({\tt G})$  generated by dominating connected subgraphs of $\tG$.

To facilitate the calculations in the remaining examples, we will use some basic notions of algebraic Morse theory, as introduced by Forman~\cite{forman1998morse}; for an overview, the reader is referred to~\cite[Section~11]{Kozlov}. Roughly speaking, algebraic Morse theory gives a convenient way to reduce a (co)chain complex by eliminating acyclic summands via changes of bases.
 
The main construction we will use goes as follows. Let $\bK$ be a field; consider a finitely generated chain complex of  $\bK$-vector spaces, say $(C_*,\partial_*)$, and a basis $B_i = \{ b_{i}^{j}\}_{j= 0,..., k_i}$ of $C_{i}$ as a $\bK$-vector space, for each $i$. With respect to these bases, the differential can be expressed as
\[ \partial(b_i^j) = \sum_{h} \langle \partial b_i^j,b_{i+1}^{h}\rangle\  b_{i+1}^{h} \ ,\]
for some coefficients $ \langle \partial b_i^j,b_{i+1}^{h}\rangle\in \bK$. One can now construct a directed graph {\tt C} with vertices $V({\tt C})\coloneqq \bigcup_{i} B_{i}$, and directed edges~$(b_{i}^{k},b_{j}^{h})\in E(\tC)$ if and only if $ \langle \partial b_i^j,b_{i+1}^{h}\rangle \neq 0$.

\begin{defn}[{\cite[Section~3]{chari2000discrete}}]\label{def:matching}
	A \emph{matching}~$M$ on a directed graph {\tt C} is  a subset of pairwise disjoint edges of {\tt C}. A matching is called \emph{acyclic} if the graph obtained from \tC~by inverting the orientations of the edges in $M$ has no directed cycles.
\end{defn}

For a chain complex $(C_*,\partial_*)$, a (acyclic) matching on $(C_*,\partial_*)$ is defined as a (acyclic) matching on the associated graph \tC.
The main result in algebraic Morse theory  (\cite[Section~8]{forman1998morse}, see also \cite{chari2000discrete}) is that, given an acyclic matching $M$ on {\tt C}, the complex $(C_*,\partial_*)$ is quasi-isomorphic to a complex~$(C_*^{M},\partial^M_*)$. Here $C_i^{M}$ is generated by all the $b_i^j$'s  that are not incident to the edges in~$M$; the generators that are not paired by $M$ are said to be \emph{critical generators} with respect to $M$. 
If $b_i^j$ and $b_{i+1}^k$ are critical generators, then $\langle \partial_*^M(b_i^j), b_{i+1}^k \rangle$ is determined by a (weighted) count of certain oriented paths called \emph{zig-zags}\footnote{These directed loops also appear in the literature as V-paths, or alternating paths. In these paths, the edges corresponding to $M$ are reversed.} in the graph $\tC$ joining  $b_i^j$ and $b_{i+1}^k$ -- see \cite[Definitions~11.1 and 11.23, Equation~(11.7)]{Kozlov}. This might make cumbersome the definition of the differential on the new complex; however, in the case at hand it is immediate to determine $\partial^M$.

\begin{rem}\label{rmk:induced diff koz}
Necessary conditions to have non-trivial zig-zags are the existence of either:
\begin{itemize}[label = {\rm (\arabic*)}]
\item[{\rm ($\diamondsuit$)}]\label{rem:pointa} a critical generator $b$, and a matched generator $a$, such that $\langle \partial(b), a \rangle \neq 0$;
\item[{\rm ($\heartsuit$)}]\label{rem:pointb} a critical generator $b'$, and a matched generator $a'$, such that $\langle \partial(a'), b' \rangle \neq 0$.
\end{itemize}
Moreover, if either the critical generators or the matched generators form a sub-complex of $(C_*,\partial_*)$, then at least one between {\rm ($\diamondsuit$)} or {\rm ($\heartsuit$)}, respectively, are violated. 
It follows that in both cases the differential~$\partial^{M}_*$ is the  restriction of $\partial_*$ to the span of the critical generators of $M$.
\end{rem}

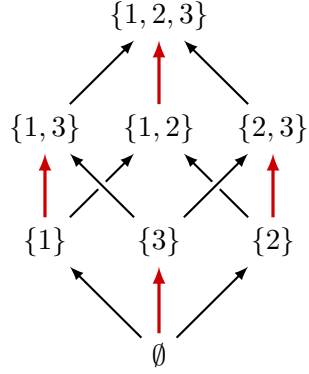
\begin{figure}[ht]
    \centering
    \begin{tikzpicture}[scale =1.5]
         \node (a) at (-1,0) {$\{1\}$};
         \node (b) at (0,1) {$\{1,2\}$};
         \node (c) at (0,-1){$\emptyset$};
         \node (d) at (1,0) {$\{2\}$};
\draw[thick, -latex] (a)-- (b);
\draw[thick, -latex] (c)-- (a);
\draw[thick, -latex] (d)-- (b);
\draw[thick, -latex] (c)-- (d);
         \node (a1) at (-1,1) {$\{1,3\}$};
         \node (b1) at (0,2) {$\{1,2,3\}$};
         \node (c1) at (0,0){$\{3\}$};
         \node (d1) at (1,1) {$\{2,3\}$};
\draw[line width = 4, white] (a1)-- (b1);
\draw[line width = 4, white] (a1)-- (c1);
\draw[line width = 4, white] (d1)-- (b1);
\draw[line width = 4, white] (d1)-- (c1);
\draw[thick, -latex] (a1)-- (b1);
\draw[thick, latex-] (a1)-- (c1);
\draw[thick, -latex] (d1)-- (b1);
\draw[thick, latex-] (d1)-- (c1);

\draw[ very thick, bunired, -latex] (d)-- (d1);
\draw[very thick, bunired, -latex] (a)-- (a1);
\draw[very thick, bunired, -latex] (b)-- (b1);
\draw[very thick, bunired, -latex] (c)-- (c1);
\end{tikzpicture}
    \caption{An acyclic matching (in red), involving all vertices, on $B(3) \cong C_*(\Delta^3)$. Inverting the orientation of the red matched edges does not introduce any directed cycle.}
    \label{fig:acy match bool}
\end{figure}

\begin{rem}\label{rmk:booleanacyclic}
Let $\tC_n$ be the face poset of a $n$-simplex $\Delta^n$. Then $\tC_n$ can be naturally identified with the Boolean poset $B(n)$. This is the graph associated to $C_{*}(\Delta^n)$. Then, an acyclic matching on $\tC_n$, involving all vertices, does always exist -- see, Figure~\ref{fig:acy match bool}, for an example. These perfect acyclic matchings correspond to a sequence of elementary collapses of the  $n$-simplex to a point.
\end{rem}

In order to facilitate the following discussions, we will make extensive use of this next result.
\begin{lem}[Technical Lemma]\label{lem:tech lemma}
Let $(C_*,\partial_*)$ be a finitely generated chain complex over a field~$\bK$. 
For each $i$, denote by $B_i = \{ b_{i}^j \}_{j=1,...,k_i}$ a basis for $C_i$ as a  $\bK$-vector space. Assume there exists a matching $M$ on $(C_*, \partial_*)$ partitioned into acyclic sub-matchings $M_1, ..., M_k$, and a function $\varphi \colon \bigcup B_i  \to \bN$ such that:
\begin{enumerate}[label = {\rm (\arabic*)}]
\item if $(a,b) \in M_s$ for some $1 \le s \le k$, then $\varphi(a) = \varphi(b)$;
\item if  $(a,b) \notin M_s$, and~$\langle \partial a,b\rangle \neq 0$ then~$\varphi(b) > \varphi(a)$.
\end{enumerate}
Then $M$ is acyclic.
\end{lem}

\begin{proof}
Assume $M$ supports a directed cycle $\gamma$; then by \cite[Lemma~2]{celoria2020filtered} the vertices of $\gamma$ must consist of a sequence of generators   $$b^{j_1}_{i},b^{r_1}_{i+1},..., b^{j_s}_{i},b^{r_s}_{i+1},$$
such that: $(b^{r_h}_{i+1},b^{j_h}_{i})\in M$, $\langle \partial b^{r_h}_{i+1}, b^{j_{h+1}}_{i}\rangle \neq 0$ for all $h = 1,\ldots,s$, and $\langle \partial b^{r_s}_{i+1}, b^{j_1}_{i}\rangle \neq 0$.
Note that $\varphi(b^{r_u}_{i+1}) \geq \varphi (b^{j_u}_{i}) \geq \varphi ( b^{r_{u+1}}_{i+1})$ for all $u =1 ,...,s$ (where $u$ is considered modulo $s$). This implies that
\begin{equation}
\varphi(b^{j_{h+1}}_{i}) \geq \dots \geq \varphi(b^{r_s}_{i+1}) \geq  \varphi(b^{j_1}_{i}) \geq \cdots \geq \varphi(b^{r_h}_{i+1})  .
\label{eq:long ineq}
\end{equation}

We have two cases; either $(b^{r_h}_{i+1},b^{j_h}_{i})\in M_t$  for all $h$ and some $t\in \{1,\dots,k\}$, or there is at least one $h$ such that  $(b^{r_h}_{i+1},b^{j_h}_{i})\in M_p$ and~$(b^{r_{h+1}}_{i+1},b^{j_{h+1}}_{i})\in M_q$ with $p\neq q$. The former case is absurd since each $M_t$ was assumed to be acyclic; for the latter, observe that $\varphi(b^{r_h}_{i+1}) > \varphi(b^{j_{h+1}}_{i})$ by hypothesis (since $M_p\cap M_q = \emptyset$), which contradicts \eqref{eq:long ineq}.
\end{proof}

\begin{figure}[ht]
\centering
\includegraphics[width=8cm]{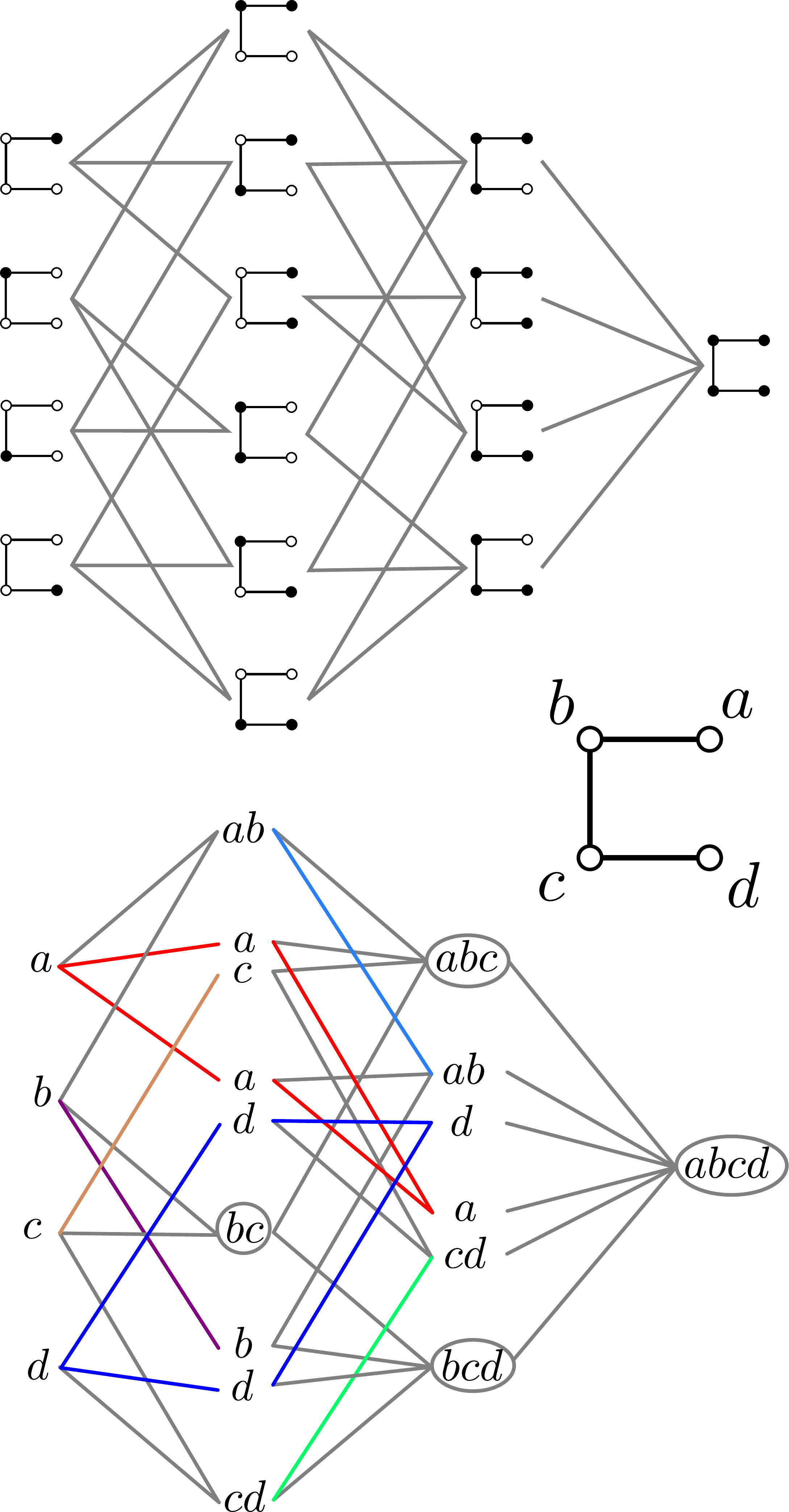}
\caption{In the top part of the figure, the Boolean poset associated to colourings of the length $3$ path (shown on the right). Below, the subdivision into connected components induced by the bi-colouring; here coloured edges provide the Boolean posets associated with each connected component. Circled elements represent (unpaired) connected dominating sets.}
\label{fig:length3path}
\end{figure}

We can regard each $x\in\pi_0(\tG_\e)$, for $\e \in \cl(\tG)$, as a connected subgraph $\tH_x$ of $\tG$  induced by its vertices. Let ${\rm D}\mathbb{H}({\tt G})$ be the sub-complex of ${\rm C}\bH(\tG)$ spanned by those $x\in\pi_0(\tG_\e)$ for which~$\tH_x$ is a dominating  connected subgraph of~$\tG$ (\emph{i.e.}~induced by a dominating set). Note that this is indeed a sub-complex, since adding any vertex to a connected dominating set induces a connected dominating subgraph. We will sometimes identify the generators of ${\rm D}\bH$ with the corresponding dominating sets, rather than with the associated dominating subgraphs. 

We can now prove a categorified version of Proposition~\ref{prop:russo};  there we showed that the Euler characteristic of $\bH$ coincides with an alternating sum of the cardinalities of connected dominating sets. Here, we prove that a similar statement holds at the level of the chain complexes. 
 
\begin{thm}\label{thm:retraction}
The chain complexes ${\rm C}\mathbb{H}({\tt G})$ and ${\rm D}\mathbb{H}({\tt G})$ are quasi-isomorphic.
\end{thm}

\begin{proof}
The bold chain group ${\rm C}\bH(\tG)$ is spanned by all $x\in \pi_0(\tG_{\e})$ with $\e\in \cl(\tG)$. 
For each generator $x\in \pi_0(\tG_{\e})$ define the quantity
\[ \varphi (x) =  |V(\tH_x)| \ .\]
Note that if $\langle \partial x , y \rangle \neq 0$, then $\varphi(x) \geq \varphi(y)$; in this case the equality is achieved if and only if $\tH_x = \tH_y$.
Now, we can consider the partition of the generators in ${\rm C}\bH_*(\tG)$ induced by the equivalence relation 
\[ x \sim x' \iff \tH_x = \tH_{x'}\ .\]
The function $\varphi$ is constant on each equivalence class.
Moreover, there is a natural identification of each such class with a Boolean poset; this is given by the correspondence
\[  [x] \ni x'\longleftrightarrow \e'\in \cl(\tH_x)\ , \]
where $x'\in\pi_0(\tG_{\e'})$. This poset is non trivial as long as $\tH_x$ is not dominating.  

We can pair up the generators in each class with an acyclic matching on the Boolean poset they span (following Remark~\ref{rmk:booleanacyclic}). 
The matching on each equivalence class is, by definition, acyclic. The function $\varphi$ is constant on each of these matchings, and it increases  if $\langle \partial x , y \rangle \neq 0$ and $[x]\neq [y]$. It follows from Lemma~\ref{lem:tech lemma} that the union of these matching is an acyclic matching.

Note that the critical generators with respect to these matchings are exactly those $x$ such that~$\tH_x$ is dominating. Since $\tH_x$ is connected and induced by its vertices, and $\cl(\tG;\tH_x)$ consists of a single colour, we can identify these generators with the corresponding connected dominating sets of $\tG$.

To conclude, the sequence of retractions induced by the acyclic matching provided above, together with Remark~\ref{rmk:induced diff koz}, give a quasi-isomorphism between  $ {\rm C}\mathbb{H}({\tt G})$ and ${\rm D}\bH(\tG)$.
\end{proof}

In particular, the homology $\bH (\tG)$ always admits a set of generators consisting of (formal sums of) connected dominating subgraphs in $\tG$.
This last result implies a few immediate corollaries.
\begin{cor}\label{cor:support}
The homology $\bH_i(\tG)$ can be non-trivial only for $\gamma^c(\tG) \leq i \leq |V(\tG)|$.
\end{cor}
\begin{proof}
By definition, there are no connected dominating sets in $\tG$ with strictly less than $\gamma^c(\tG)$ vertices.
\end{proof}

\begin{cor}
If \tG~is a disconnected simple graph, then  $\bH({\tt G}) = 0$.
\end{cor}
\begin{proof}
As the graph is not connected, there are no connected dominating sets.
\end{proof}

We conclude this subsection by remarking that the description of $\bH$ using dominating sets can be used, in conjunction with techniques from algebraic Morse theory, to provide an alternative proof of Proposition~\ref{prop: complete graphs}.

\begin{proof}[Alternative proof of Proposition~\ref{prop: complete graphs}]\label{rem:char of Km II: the vengeance}
We  claim that if $\tG$ is not a complete graph, then $\bH_1(\tG)=0$.
By definition, ${\rm D}\bH_1(\tG)$ is spanned by all dominating vertices of $\tG$, and ${\rm D}\bH_2(\tG)$ contains all edges in~$E(\tG)$ such that at least one of their endpoints is a dominating vertex. 
If there are no dominating vertices, then ${\rm D}\bH_1(\tG) = 0$, and we are done.
Otherwise, $\tG$ must have a non-dominating vertex~$v$. We can consider the matching given by $(w,\{ w, v\})$, for each dominating vertex $w$ in $\tG$.
The matching so constructed is acyclic since, in the graph associated to ${\rm D}\bH(\tG)$, $(w,\{ w, v\})$ is the only edge with target $\{ w, v\}$. 
As there are no critical generators in degree $1$, the claim follows from \cite[Theorem~11.24]{Kozlov}.
\end{proof}

\section{Applications and computations}\label{sec: appl and comps}

In this final section we collect some consequences stemming from Theorem~\ref{thm:retraction}, as well as some computations (over a field).

\subsection{Applications}
In this first subsection we re-obtain the full computation of the bold homology for complete graphs, and present some vanishing results stemming from Theorem~\ref{thm:retraction}. In particular,  we can also easily deduce the next result (\emph{cf.}~\cite[Conjecture~8.2]{uberhomology}); remarkably, this is also a consequence of the more general Proposition~\ref{prop:leaf} proved below.
\begin{thm}
Let $\tT$ be a connected tree on $n\geq 3$ vertices. Then $\bH(\tT) = 0.$
\end{thm}
\begin{proof}
The connected dominating sets in a connected tree $\tT$ are easily seen to be those obtained by discarding any number of univalent vertices (possibly all, since $n\geq 3$). If $\tT$ has $m$ leaves, then~${\rm D}\bH(\tT)$ is isomorphic to the simplicial chain complex $C(\Delta^m)$ of the $m$-simplex $\Delta^m$. Therefore, $\bH(\tT)$ vanishes --\emph{cf.}~Remark~\ref{rmk:booleanacyclic}.
\end{proof}

\begin{figure}[ht]
    \centering
    \begin{tikzpicture}[baseline=(current bounding box.center)]
		\tikzstyle{point}=[circle,thick,draw=black,fill=black,inner sep=0pt,minimum width=2pt,minimum height=2pt]
		\tikzstyle{arc}=[shorten >= 8pt,shorten <= 8pt,->, thick]
		
		\node at (-1,-2) {$\tt{G}'$};
		\draw[dotted] (-1,-2)  circle (.7);
		\node[above] (v0) at (0,-2) {$v$};
		\draw[fill] (0,-2)  circle (.05);
		\node[above] (v1) at (1.5,-2) {$w$};
		\draw[fill] (1.5,-2)  circle (.05);
		\node at (-.4,-1.9) {\tiny{$\vdots$}};
		
		\draw[thick, bunired] (-0.35,-2.25) -- (-0.15,-2.05);
		\draw[thick, bunired] (-0.35,-1.75) -- (-0.15,-1.95);
		\draw[thick, bunired] (0.15,-2) -- (1.35,-2);
	\end{tikzpicture}
    \caption{A leaf in the graph $\tG$ induces the splitting shown above.}
    \label{fig:leaf}
\end{figure}
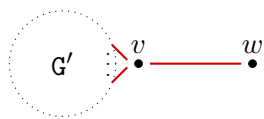

More generally, we obtain a vanishing result for graphs with a leaf (that is a univalent vertex, \emph{cf.}~Figure~\ref{fig:leaf}):

\begin{prop}\label{prop:leaf}
Let $\tG$ be a simple and connected graph on $n\geq 3$ vertices. If $\tG$ contains a leaf, then $\mathbb{H}(\tG) = 0$.
\end{prop}
\begin{proof}

Since $\tG$ has a leaf, it must contain a vertex $v$ and a univalent vertex $w$, as in Figure~\ref{fig:leaf}.  
Note that any connected dominating set of $\tG$ must contain the vertex $v$.
In this case, a matching can be described explicitly: consider a connected and dominating set $D$ of $\tG' \coloneqq \tG\setminus \{ w \}$ containing $v$; this is also a dominating set of $\tG$. We can then  pair $D$ with the dominating set $D\cup \{ w\}$. 

Since $n\geq 3$, these are all the possible dominating sets of $\tG$. This matching is easily verified to be acyclic, and the statement follows.
\end{proof}

\begin{figure}[ht]
		\centering
		\newdimen\R
		\R=2.0cm
		\begin{subfigure}[b]{0.3\textwidth}
		 \centering
		\begin{tikzpicture}
			\draw[xshift=5.0\R, fill] (270:\R) circle(.05)  node[below] {$v_n$};
			\draw[xshift=5.0\R,fill] (225:\R) circle(.05)  node[below left]   {$v_1$};
			\draw[xshift=5.0\R,fill] (180:\R) circle(.05)  node[left] {$v_2$};
			\draw[xshift=5.0\R,fill] (135:\R) circle(.05)  node[above left] {$v_3$};
			\draw[xshift=5.0\R, fill] (90:\R) circle(.05)  node[above] {$v_4$};
			\draw[xshift=5.0\R,fill] (45:\R) circle(.05)  node[above right] {$v_5$};
			\draw[xshift=5.0\R,fill] (0:\R) circle(.05)  node[right] {$v_6$};
			\draw[xshift=5.0\R,fill] (315:\R) circle(.05)  node[below right] {$v_{n-1}$};
			
			\node[xshift=5.0\R] (p) at (0,0) { };
			\node[xshift=5.0\R] (v0) at (270:\R) { };
			\node[xshift=5.0\R] (v1) at (225:\R) { };
			\node[xshift=5.0\R] (v2) at (180:\R) { };
			\node[xshift=5.0\R] (v3) at (135:\R) { };
			\node[xshift=5.0\R] (v4) at (90:\R) { };
			\node[xshift=5.0\R] (v5) at (45:\R) { };
			\node[xshift=5.0\R] (v6) at (0:\R) { };
			\node[xshift=5.0\R] (vn) at (315:\R) { };
			
			\draw[thick, red] (v0)--(v1);
			\draw[thick, red] (v1)--(v2);
			\draw[thick, red] (v2)--(v3);
			\draw[thick, red] (v3)--(v4);
			\draw[thick, red] (v4)--(v5);
			\draw[thick, red] (v5)--(v6);
			\draw[thick, red] (vn)--(v0);

			\draw[xshift=5.0\R, fill] (292.5:\R) node[below right] {};
			\draw[xshift=5.0\R,fill] (247.5:\R) node[below left] {};
			\draw[xshift=5.0\R,fill] (202.5:\R)   node[left] {};
			\draw[xshift=5.0\R,fill] (157.5:\R)  node[above left] {};
			\draw[xshift=5.0\R, fill] (112.5:\R)   node[above] {};
			\draw[xshift=5.0\R,fill] (67.5:\R) node[above right] {};
			\draw[xshift=5.0\R,fill] (22.5:\R) node[right] {};
			\draw[xshift=4.95\R,fill] (337.5:\R)  node {$\cdot$} ;
			\draw[xshift=4.95\R,fill] (333:\R)  node {$\cdot$} ;
			\draw[xshift=4.95\R,fill] (342:\R)  node {$\cdot$} ;

			\node[xshift=5.0\R, above right] (vhat) at (0.15,0.1) {};
		\end{tikzpicture}
		\caption{}
		\end{subfigure}
		\hspace{.1\textwidth}
		\begin{subfigure}[b]{.3\textwidth}
		\begin{tikzpicture}
			\draw[xshift=5.0\R, fill] (270:\R) circle(.05)  node[below] {$v_n$};
			\draw[xshift=5.0\R,fill] (225:\R) circle(.05)  node[below left]   {$v_1$};
			\draw[xshift=5.0\R,fill] (180:\R) circle(.05)  node[left] {$v_2$};
			\draw[xshift=5.0\R,fill] (135:\R) circle(.05)  node[above left] {$v_3$};
			\draw[xshift=5.0\R, fill] (90:\R) circle(.05)  node[above] {$v_4$};
			\draw[xshift=5.0\R,fill] (45:\R) circle(.05)  node[above right] {$v_5$};
			\draw[xshift=5.0\R,fill] (0:\R) circle(.05)  node[right] {$v_6$};
			\draw[xshift=5.0\R,fill] (315:\R) circle(.05)  node[below right] {$v_{n-1}$};
			
			\node[xshift=5.0\R] (p) at (0,0) { };
			\node[xshift=5.0\R] (v0) at (270:\R) { };
			\node[xshift=5.0\R] (v1) at (225:\R) { };
			\node[xshift=5.0\R] (v2) at (180:\R) { };
			\node[xshift=5.0\R] (v3) at (135:\R) { };
			\node[xshift=5.0\R] (v4) at (90:\R) { };
			\node[xshift=5.0\R] (v5) at (45:\R) { };
			\node[xshift=5.0\R] (v6) at (0:\R) { };
			\node[xshift=5.0\R] (vn) at (315:\R) { };
			
			\draw[thick, red] (v0)--(v1);
			\draw[thick, red] (v1)--(v2);
			\draw[thick, red] (v2)--(v3);
			\draw[thick, red] (v3)--(v4);
			\draw[thick, red] (v4)--(v5);
			\draw[thick, red] (v5)--(v6);
			\draw[thick, red] (vn)--(v0);

			\draw[thick, blue] (v0)--(p);			
			\draw[thick, blue] (vn)--(p);
			\draw[thick, blue] (v1)--(p);			
			\draw[thick, blue] (v2)--(p);
			\draw[thick, blue] (v3)--(p);			
			\draw[thick, blue] (v4)--(p);
			\draw[thick, blue] (v5)--(p);			
			\draw[thick, blue] (v6)--(p);
			
			\draw[xshift=5.0\R, fill] (292.5:\R) node[below right] {};
			\draw[xshift=5.0\R,fill] (247.5:\R) node[below left] {};
			\draw[xshift=5.0\R,fill] (202.5:\R)   node[left] {};
			\draw[xshift=5.0\R,fill] (157.5:\R)  node[above left] {};
			\draw[xshift=5.0\R, fill] (112.5:\R)   node[above] {};
			\draw[xshift=5.0\R,fill] (67.5:\R) node[above right] {};
			\draw[xshift=5.0\R,fill] (22.5:\R) node[right] {};
			\draw[xshift=4.95\R,fill] (337.5:\R)  node {$\cdot$} ;
			\draw[xshift=4.95\R,fill] (333:\R)  node {$\cdot$} ;
			\draw[xshift=4.95\R,fill] (342:\R)  node {$\cdot$} ;

			\node[xshift=5.0\R, above right] (vhat) at (0.15,0.1) {};
			\draw[fill, white] (vhat) circle (.15);
			
			\draw[xshift=5.0\R,fill] (0,0) circle(.05)  node[above right] {$\hat{v}$};
		\end{tikzpicture}\caption{}
		\end{subfigure}
		\caption{The cycle graph~$\tC_n$ and the wheel graph $\tW_{n+1} = \mathrm{Cone}(\tC_n)$.} 
		\label{fig:polycone}
	\end{figure}
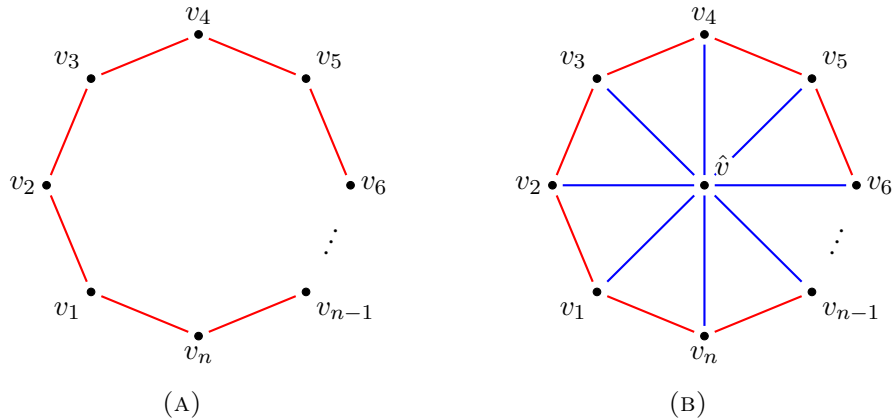

Graphs with leaves are not the only graphs for which the bold homology vanishes. In order to see this, we first need to investigate the behaviour of bold homology with respect to the cone operation. For a graph $\tG$, we consider its \emph{graph cone} ${\rm Cone}(\tG)$. This is the graph obtained by adding one extra vertex $\hat{v}$ to the vertices of $\tG$, and one edge between $\hat{v}$ and each  $v\in V(\tG)$. As an example, consider the graph in Figure~\ref{fig:polycone}, or in Figure~\ref{fig:gem}.

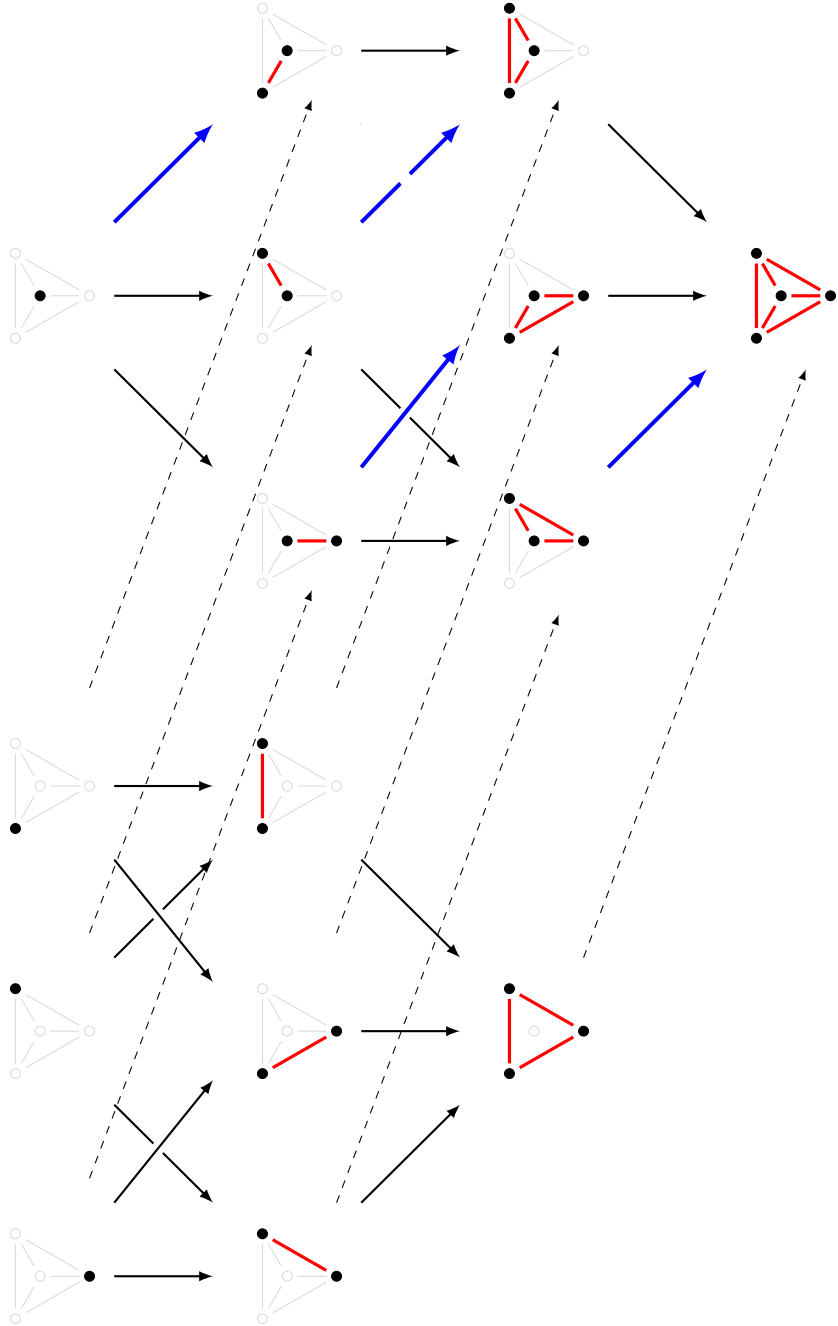
\begin{figure}[ht]
    \centering
        \begin{tikzpicture}[scale =.65,thin]

\begin{scope}[shift = {+(-5,-15)}]

    \begin{scope}[shift = {+(5,0)}]
        \node[] (a) at (0,0){};
        \draw[gray, opacity = .3] (a) circle (.1);
        \foreach \a in {1,...,3}
         {
            \node[] (u\a) at ({\a*120}:1) {};
            \draw[gray, opacity = .3] (a) -- (u\a);
            \draw[gray, opacity = .3] (u\a) circle (.1);
         } 
         \draw[fill] (u1) circle (.1);
         \draw[gray, opacity = .3] (u1) -- (u2) -- (u3) -- (u1);
    \end{scope}
    
    \begin{scope}[shift = {+(5,5)}]
        \node[] (a) at (0,0){};
        \draw[gray, opacity = .3] (a) circle (.1);
        \foreach \a in {1,...,3}
         {
            \node[] (u\a) at ({\a*120}:1) {};
            \draw[gray, opacity = .3] (a) -- (u\a);
            \draw[gray, opacity = .3] (u\a) circle (.1);
         }
         \draw[fill] (u2) circle (.1);
         \draw[gray, opacity = .3] (u1) -- (u2) -- (u3) -- (u1);
    \end{scope}

    \begin{scope}[shift = {+(5,-5)}]
        \node[] (a) at (0,0){};
        \draw[gray, opacity = .3] (a) circle (.1);
        \foreach \a in {1,...,3}
         {
            \node[] (u\a) at ({\a*120}:1) {};
            \draw[gray, opacity = .3] (a) -- (u\a);
            \draw[gray, opacity = .3] (u\a) circle (.1);
         }
         \draw[fill] (u3) circle (.1);
         \draw[gray, opacity = .3] (u1) -- (u2) -- (u3) -- (u1);
    \end{scope}
    
     \begin{scope}[shift = {+(10,0)}]
        \node[] (a) at (0,0){};
        \draw[gray, opacity = .3] (a) circle (.1);
        \foreach \a in {1,...,3}
         {
            \node[] (u\a) at ({\a*120}:1) {};
            \draw[gray, opacity = .3] (a) -- (u\a);
            \draw[gray, opacity = .3] (u\a) circle (.1);
         }
         \draw[fill] (u3) circle (.1);
         \draw[fill] (u2) circle (.1);
         \draw[gray, opacity = .3] (u1) -- (u2) -- (u3) -- (u1);
         \draw[very thick, red] (u2) -- (u3);
    \end{scope}
    
    \begin{scope}[shift = {+(10,5)}]
        \node[] (a) at (0,0){};
        \draw[gray, opacity = .3] (a) circle (.1);
        \foreach \a in {1,...,3}
         {
            \node[] (u\a) at ({\a*120}:1) {};
            \draw[gray, opacity = .3] (a) -- (u\a);
            \draw[gray, opacity = .3] (u\a) circle (.1);
         }
         \draw[fill] (u1) circle (.1);
         \draw[fill] (u2) circle (.1);
         \draw[gray, opacity = .3] (u1) -- (u2) -- (u3) -- (u1);
         \draw[very thick, red] (u2) -- (u1);
    \end{scope}

    \begin{scope}[shift = {+(10,-5)}]
        \node[] (a) at (0,0){};
        \draw[gray, opacity = .3] (a) circle (.1);
        \foreach \a in {1,...,3}
         {
            \node[] (u\a) at ({\a*120}:1) {};
            \draw[gray, opacity = .3] (a) -- (u\a);
            \draw[gray, opacity = .3] (u\a) circle (.1);
         } 
          \draw[fill] (u3) circle (.1);
         \draw[fill] (u1) circle (.1);
         \draw[gray, opacity = .3] (u1) -- (u2) -- (u3) -- (u1);
         \draw[very thick, red] (u1) -- (u3);
    \end{scope}
    \begin{scope}[shift = {+(15,0)}]
        \node[] (a) at (0,0){};
        \draw[gray, opacity = .3] (a) circle (.1);
        \foreach \a in {1,...,3}
         {
            \node[] (u\a) at ({\a*120}:1) {};
            \draw[fill] (u\a) circle (.1);
         } 
         \draw[very thick, red] (u1) -- (u2) -- (u3) -- (u1);
    \end{scope}

    \draw[thick, -latex] (6.5,1.5) -- (8.5,3.5) ;
    \draw[thick, -latex] (6.5,-1.5) -- (8.5,-3.5) ;
    
    \draw[line width =5, white ] (6.5,3.5) -- (8.5,1) ;
    \draw[thick, -latex] (6.5,3.5) -- (8.5,1) ;
    \draw[line width =5, white ] (6.5,-3.5) -- (8.5,-1) ;
    \draw[thick, -latex] (6.5,-3.5) -- (8.5,-1) ;
    
    \draw[thick, -latex] (6.5,5) -- (8.5,5) ;
    \draw[thick, -latex] (6.5,-5) -- (8.5,-5) ;
    
    \draw[thick, latex-] (13.5,1.5) -- (11.5,3.5) ;
    \draw[thick, latex-] (13.5,-1.5) -- (11.5,-3.5) ;
    \draw[thick, -latex] (11.5,0) -- (13.5,0) ;
\end{scope}

\draw[thin, dashed, -latex] (11,-13.5) -- (15.5,-1.5) ;

\draw[thin, dashed, -latex] (6,-18.5) -- (10.5,-6.5) ;
\draw[thin, dashed, -latex] (6,-13) -- (10.5,-1) ;
\draw[thin, dashed, -latex] (6,-8) -- (10.5,4) ;

\begin{scope}[shift = {+(-5,0)}]
\draw[thin,dashed, -latex] (6,-18) -- (10.5,-6) ;
\draw[thin,dashed, -latex] (6,-13) -- (10.5,-1) ;
\draw[thin,dashed, -latex] (6,-8) -- (10.5,4) ;
\end{scope}

    \begin{scope}[shift = {+(0,0)}]
        \node[] (a) at (0,0){};
        \draw[fill] (a) circle (.1);
        \foreach \a in {1,...,3}
         {
            \node[] (u\a) at ({\a*120}:1) {};
            \draw[gray, opacity = .3] (a) -- (u\a);
            \draw[gray, opacity = .3] (u\a) circle (.1);
         }
         \draw[gray, opacity = .3] (u1) -- (u2) -- (u3) -- (u1);
    \end{scope}
    
    \begin{scope}[shift = {+(5,0)}]
        \node[] (a) at (0,0){};
        \draw[fill] (a) circle (.1);
        \foreach \a in {1,...,3}%
         {
            \node[] (u\a) at ({\a*120}:1) {};
            \draw[gray, opacity = .3] (a) -- (u\a);
            \draw[gray, opacity = .3] (u\a) circle (.1);
         } 
         \draw[fill] (u1) circle (.1);
         \draw[gray, opacity = .3] (u1) -- (u2) -- (u3) -- (u1);
         \draw[very thick, red] (u1) -- (a);
    \end{scope}
    
    \begin{scope}[shift = {+(5,5)}]
        \node[] (a) at (0,0){};
        \draw[fill] (a) circle (.1);
        \foreach \a in {1,...,3}
         {
            \node[] (u\a) at ({\a*120}:1) {};
            \draw[gray, opacity = .3] (a) -- (u\a);
            \draw[gray, opacity = .3] (u\a) circle (.1);
         } 
         \draw[fill] (u2) circle (.1);
         \draw[gray, opacity = .3] (u1) -- (u2) -- (u3) -- (u1);
         \draw[very thick, red] (u2) -- (a);
    \end{scope}

    \begin{scope}[shift = {+(5,-5)}]
        \node[] (a) at (0,0){};
        \draw[fill] (a) circle (.1);
        \foreach \a in {1,...,3}
         {
            \node[] (u\a) at ({\a*120}:1) {};
            \draw[gray, opacity = .3] (a) -- (u\a);
            \draw[gray, opacity = .3] (u\a) circle (.1);
         } 
         \draw[fill] (u3) circle (.1);
         \draw[gray, opacity = .3] (u1) -- (u2) -- (u3) -- (u1);
         \draw[very thick, red] (u3) -- (a);
    \end{scope}
    
     \begin{scope}[shift = {+(10,0)}]
        \node[] (a) at (0,0){};
        \draw[fill] (a) circle (.1);
        \foreach \a in {1,...,3}
         {
            \node[] (u\a) at ({\a*120}:1) {};
            \draw[gray, opacity = .3] (a) -- (u\a);
            \draw[gray, opacity = .3] (u\a) circle (.1);
         } 
         \draw[fill] (u3) circle (.1);
         \draw[fill] (u2) circle (.1);
         \draw[gray, opacity = .3] (u1) -- (u2) -- (u3) -- (u1);
         \draw[very thick, red] (u2) -- (u3) -- (a) -- (u2);
    \end{scope}
    
    \begin{scope}[shift = {+(10,5)}]
        \node[] (a) at (0,0){};
        \draw[fill] (a) circle (.1);
        \foreach \a in {1,...,3}
         {
            \node[] (u\a) at ({\a*120}:1) {};
            \draw[gray, opacity = .3] (a) -- (u\a);
            \draw[gray, opacity = .3] (u\a) circle (.1);
         } 
         \draw[fill] (u1) circle (.1);
         \draw[fill] (u2) circle (.1);
         \draw[gray, opacity = .3] (u1) -- (u2) -- (u3) -- (u1);
         \draw[very thick, red] (u2) -- (u1) -- (a) -- (u2);
    \end{scope}

    \begin{scope}[shift = {+(10,-5)}]
        \node[] (a) at (0,0){};
        \draw[fill] (a) circle (.1);
        \foreach \a in {1,...,3}
         {
            \node[] (u\a) at ({\a*120}:1) {};
            \draw[gray, opacity = .3] (a) -- (u\a);
            \draw[gray, opacity = .3] (u\a) circle (.1);
         } 
          \draw[fill] (u3) circle (.1);
         \draw[fill] (u1) circle (.1);
         \draw[gray, opacity = .3] (u1) -- (u2) -- (u3) -- (u1);
         \draw[very thick, red] (u1) -- (u3) -- (a) -- (u1);
    \end{scope}
    \begin{scope}[shift = {+(15,0)}]
        \node[] (a) at (0,0){};
        \draw[fill] (a) circle (.1);
        \foreach \a in {1,...,3}
         {
            \node[] (u\a) at ({\a*120}:1) {};
            \draw[very thick, red] (a) -- (u\a);
            \draw[fill] (u\a) circle (.1);
         } 
         \draw[very thick, red] (u1) -- (u2) -- (u3) -- (u1);
    \end{scope}
    \draw[line width = 1.5 , blue, -latex] (1.5,1.5) -- (3.5,3.5) ;
    \draw[thick, -latex] (1.5,-1.5) -- (3.5,-3.5) ;
    \draw[thick, -latex] (1.5,0) -- (3.5,0) ;
       \draw[thick, -latex] (6.5,3.5) -- (8.5,1) ;
    \draw[line width =5, white ] (6.5,1.5) -- (8.5,3.5) ;
    \draw[line width = 1.5 , blue, -latex] (6.5,1.5) -- (8.5,3.5) ;
    \draw[thick, -latex] (6.5,-1.5) -- (8.5,-3.5) ;
    
    \draw[line width =5, white ] (6.5,3.5) -- (8.5,1) ;
 
    \draw[line width =5, white ] (6.5,-3.5) -- (8.5,-1) ;
    \draw[line width = 1.5 , blue, -latex] (6.5,-3.5) -- (8.5,-1) ;
    
    \draw[thick, -latex] (6.5,5) -- (8.5,5) ;
    \draw[thick, -latex] (6.5,-5) -- (8.5,-5) ;
    
    \draw[thick, latex-] (13.5,1.5) -- (11.5,3.5) ;
    \draw[line width = 1.5 , blue, latex-] (13.5,-1.5) -- (11.5,-3.5) ;
    \draw[thick, -latex] (11.5,0) -- (13.5,0) ;

\end{tikzpicture}

    \caption{The chain complex for the cone of $\tC_3$; blue edges represent an acyclic matching on the Boolean component of the complex. Dashed arrows indicate the components of the differential joining a critical generator with a matched generator.}
    \label{fig:cono}
\end{figure}

\begin{prop}\label{prop:cone}
If $\tG$ is a simple and connected graph, then $\bH({\rm Cone}(\tG)) \cong \bH(\tG)$.
\end{prop}
\begin{proof}
We define an acyclic matching $M$ on ${\rm C}\bH({\rm Cone}(\tG))$. Dominating sets in ${\rm Cone}(\tG)$ can be divided in two classes: those which contain $\hat{v}$, and those which do not (see, for an example, Figure~\ref{fig:cono}). The latter kind can be identified with the dominating sets of $\tG$ (seen as a subgraph of ${\rm Cone}(\tG)$).  Any set of vertices containing $\hat{v}$ is connected and dominating. 

The subgraphs containing $\hat{v}$ form a sub-complex, which is isomorphic to  the simplicial chain complex $C_*(\Delta^{\vert V(\tG)\vert})$ of the standard $\vert V(\tG)\vert$-simplex $\Delta^{\vert V(\tG)\vert}$; we can now consider one of the acyclic matchings whose existence is guaranteed by Remark~\ref{rmk:booleanacyclic}.

The complex $({\rm C}\bH({\rm Cone}(\tG)))^M$ induced by the critical generators with respect to $M$ (see Definition~\ref{def:matching} and subsequent lines)
is a quotient complex of ${\rm C}\bH({\rm Cone}(\tG))$ by the sub-complex spanned by the matched generators. By identifying the corresponding connected dominating sets, we can now identify the complex $({\rm C}\bH({\rm Cone}(\tG)))^M$ with~${\rm C}\bH(\tG)$,
and conclude using Remark~\ref{rmk:induced diff koz}.
\end{proof}

As a straightforward consequence of Proposition~\ref{prop:cone}, we can easily deduce the full computation of the bold homology of complete graphs. Note that this computation was already carried out in \cite[Section~8]{uberhomology} using different techniques.
\begin{cor}
Let $\tK_n$ be the complete graph on $n$ vertices; then $\bH_*(\tK_n) \cong\bH_*(\tK_1)$. In particular, $\bH_*(\tK_n)$ is of rank one in degree one, and trivial otherwise.
\end{cor}
\begin{proof}
The result follows immediately from Proposition~\ref{prop:cone}, after noting that $\tK_n$ is the result of iterating $n$ times the cone graph construction on $\tK_0$.
\end{proof}

It follows from Proposition~\ref{prop:cone} that a necessary condition for a class of graphs to be detected by $\bH$ is to be closed under graph coning. 

\begin{rem}
The converse of Proposition~\ref{prop:leaf} is false;
by Proposition~\ref{prop:cone}, $\mathbb{H}^0(\tG) = 0$ for $\tG$ the (leafless) Gem graph shown in Figure~\ref{fig:gem} below. A further example, which is not a graph cone is the Durer graph (also known as the generalised Petersen graph~$(6,2)$); we computed its homology using the program \cite{githububer}, showing that it is trivial. 
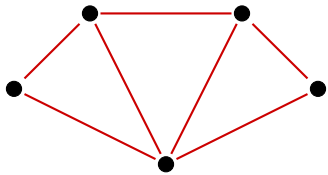
\begin{figure}[ht]
    \centering
    \begin{tikzpicture}
        \node   (v1) at (-2,-1) {};
        \node   (v2) at (-1,0) {};
        \node   (v3) at (1,0) {};
        \node   (v4) at (2,-1) {};
        \node   (v5) at (0,-2) {};
        
        \foreach \x in {1,2,3,4}{
        \draw[fill] (v\x) circle (.1); 
        \draw[thick,bunired] (v\x) -- (v5);}
        \draw[fill] (v5) circle (.1); 
        \draw[thick,bunired] (v1) -- (v2) -- (v3) -- (v4);
    \end{tikzpicture}
    \caption{The Gem graph; this is the graph cone of $\tP_3$.}
    \label{fig:gem}
\end{figure}
\end{rem}

\subsection{Computations}

We turn to some explicit computations of bold homology groups of certain families of graphs, complementing those provided in \cite{uberhomology} -- see also Table~\ref{tab:summary_table}.

We start with the computation of the bold homology groups of polygonal graphs, proving the last point of \cite[Conjecture~8.2]{uberhomology}.

\begin{prop}
Let $\tC_n$ be the cycle graph on $n$ vertices. Then, $\bH_{*} (\tC_n) \cong  \widetilde{\mathrm{H}}_{n-1-*}(S^{1})$.
\end{prop}
\begin{proof}
We can identify $\tC_n$ with a simplicial realisation of the $1$-sphere $S^1$, by identifying each $1$-simplex with the corresponding edge of $\tC_n$. 

The only connected dominating sets of $\tC_n$ are  all subsets of $V(\tC_n)$ with either $n$, $n-1$ or $n-2$ elements.
Thus, to each connected dominating set $D\subseteq V(\tC_n)$ we can associate the (possibly empty) simplex spanned by the vertices in $V(\tC)\setminus D$. This gives a bijection between the generators in ${\rm D}\bH(\tC_n)$ and those of $\widetilde{\mathrm{C}}(S^1)$,  the simplicial chain complex corresponding to the simplicial structure defined by $\tC_n$. From this description, it is immediate to see that the induced linear map $\phi$ inverts the homological degree and shifts it by $(n-1)$. Furthermore, since the differential in ${\rm D}\bH(\tC_n)$ is induced by the inclusion, it commutes with $\phi$ (up to the choice of a sign assignment, which does not affect the isomorphism class of the complex). Therefore, $\phi$ is an isomorphism of chain complexes between
${\rm D}\bH_{*} (\tC_n)$ and $ \widetilde{\mathrm{C}}_{n - 1 - *}(S^{1})$, concluding the proof. 
\end{proof}

In the above proof, we defined an explicit isomorphism $\bH_{*} (\tC_n) \cong \widetilde{\mathrm{H}}_{n - 1 - *}(\tC_n)$. This isomorphism also gives an explicit cycle in ${\rm D}\bH (\tC_n)$ whose class generates $\bH(\tC_n)$; this cycle is given by the sum of all the $(n-2)$-paths in $\tC_n$, and it has degree $n-2$.

\begin{prop}\label{prop:completebipartite}
The homology $\bH (\tK_{n,m})$ is of rank $1$ in degree $2$, and trivial otherwise.
\end{prop}
\begin{proof}
Let $V = \{v_i\}_{i = 1,\ldots, n}$ and $ W = \{w_i\}_{i = 1,\ldots, m}$ denote the two sets of vertices in $\tK_{m,n}$. A subset of $V(\tK_{m,n})$ is connected and dominating if and only if it contains at least one element in $V$ and one in  $W$.
Consider a matching $M$ on the complex ${\rm D}\bH(\tK_{n,m})$, given by the pairs $(D,D')$ of connected dominating sets such that:
\[ D' = \begin{cases} D \cup \{ v_1 \} & \text{if } v_1\notin D, \\
D \cup \{ w_1 \} & \text{if } D\cap V = \{ v_1 \} \text{ and } w_1\notin D.\end{cases}\]
We have to show that all such pairs are disjoint.
First, note that the first entry of each pair is completely determined by the second entry; given a pair $(D,D')$ either~$D\cap V = D'\cap V = \{ v_1 \}$, which implies $D = D'\setminus \{ w_1 \}$, or $D = D'\setminus \{ v_1 \}$.
Assume that, for some pairs $(D_1, D_1')$ and $(D_2, D_2')$, we have $D_1 = D_2'$. Then we must have $v_1 \notin D_1$, and $D_2\cap V =D_2'\cap V =  \{ v_1\}$.
This contradicts the fact that $D_1$ is a connected dominating set, since otherwise we would have $D_1\cap V \neq \emptyset$.
Note that by construction the only connected dominating set which is critical for $M$ is $\{ v_1, w_1\}$.

In order to conclude, we are left to show that this matching $M$ is acyclic; to this end consider the function 
\[ \varphi(D) = \begin{cases} \# D\setminus \{ w_1\} & \text{ if } D\cap V = \{ v_1\},\\
 \# D\setminus \{ v_1 \} & \text{otherwise.}\end{cases}\]
This function is non-increasing along the edges of $M$, while it increases along the other components of the differential. The statement now follows from Lemma~\ref{lem:tech lemma}.
\end{proof}

Given a pair of rooted and  connected graphs $\tG_0$ and $\tG_1$ with at least one edge, we can construct an infinite sequence of graph $NS(\tG_0,\tG_1,k)$, indexed by $k \ge 1$ as follows: $NS(\tG_0,\tG_1,1)$ is the connected graph obtained by joining the two roots with an edge $e$. Then $NS(\tG_0,\tG_1,k)$ is just obtained by subdividing $e$ $k$ time.
\begin{prop}\label{prop:neckstretch} For any pair of simple and connected graphs $\tG_0$, $\tG_1$ and integer $k \ge 1$ we have $$\bH_{*+k} (NS(\tG_0,\tG_1,k)) \cong \bH_* (NS(\tG_0,\tG_1,1)).$$
\end{prop}
\begin{proof}
The dominating and connected sets in $NS(\tG_0,\tG_1,k)$ are clearly in bijection with those in $NS(\tG_0,\tG_1,1)$. The bijection is obtained by colouring in the new vertices obtained via the subdivisions. Since the homological degree in ${\rm D}\bH$ is the number of $1$-coloured vertices, it follows that this bijection (which clearly commutes with the differential) shifts the degree by exactly~$k$.  
\end{proof}

This last result provides us with infinite families of graphs where the bold homology stabilises, up to a degree shift. On a practical level this can be used to reduce the computations of graphs with a long isthmus.  
~\\
We conclude this section with some computations, collected in Table~\ref{tab:summary_table} below of computations. We remark that in general: the rank of $\bH$ can be grater than 1; it can be supported in more than one homological degree, and it is not completely determined by its Euler characteristic. An example showing that all these facts hold is given by $\tK_3\times \tC_4$.
\begin{table}[h]
    \centering
    \begin{tabular}{l|cr}
       Graph  & $\bH$ & \raisebox{.1em}{$\chi_{\bH}$} \\
       \hline\hline
          $\tK_n$    & $\mathbb{F}_{(1)}$  & $-1$ \\
          $\tK_{m,n}$, $m,n\geq 2$    & $\mathbb{F}_{(2)}$ & $1$ \\
          $\tC_{n}$    & $\mathbb{F}_{(n-2)}$  & $(-1)^{n}$ \\
          $\tW_{n}$    & $\mathbb{F}_{(n-3)}$  & $(-1)^{n+1}$ \\
          $\tL_n$, $n>2$    &  $(0)$  & $0$ \\
          Trees $\neq \tL_2, \tL_1$    &  $(0)$  & $0$ \\
          ${\tt Cube}(2) = \tC_4$    & $\mathbb{F}_{(2)}$  & $1$ \\
          ${\tt Cube}(3) = \tC_4 \times \tL_2$    & $\mathbb{F}^{3}_{(4)}$  & $3$ \\
          ${\tt Cube}(4)$    & $\mathbb{F}^{21}_{(8)}$  & $21$ \\
          ${\tt Cube}(5)$    & unknown  & $\pm 455$ \\
          Petersen graph  & $\mathbb{F}_{(4)}$  & $1$ \\
          $\tK_3\times \tL_2$    & $\mathbb{F}_{(2)}$  & $1$ \\
          $\tK_4\times \tL_2$    & $\mathbb{F}_{(2)}$  & $1$ \\
          $\tK_3\times \tC_4$    & $\mathbb{F}_{(5)}\oplus \mathbb{F}_{(6)}^{2}$  & $1$ \\
          $\tK_4\times \tC_4$    & $\mathbb{F}_{(5)}\oplus \mathbb{F}_{(7)}^{2}$  & $-3$ \\
          $\tK_5\times \tC_4$    & unknown  & $-1$ \\
          $\tK_6\times \tC_4$    & unknown  & $-3$ \\
          $\tK_3\times \tK_3$    & $\mathbb{F}_{(4)}^5$  & $5$ \\
          $\tC_3\times \tL_2$    & $\mathbb{F}_{(2)}$  & $1$ \\
          $\tC_5\times \tL_2$    & $\mathbb{F}_{(4)}$  & $1$ \\
          $\tC_6\times \tL_2$    & $\mathbb{F}_{(6)}$  & $1$ \\
          $\tC_7\times \tL_2$    & $\mathbb{F}_{(8)}$  & $1$ \\
    \end{tabular}
    \caption{Computations of $\bH$ and/or $\chi_{\bH}$ for some specific graphs. All computations are made with coefficients in $\bF$, and we denoted by $\bF_{(i)}^{k}$ a summand of rank~$k$ in homological degree $i$.}
    \label{tab:summary_table}
\end{table}

\subsection{Open questions}\label{sec: questions}

We list here a few open questions.

\begin{q}\label{q: kunneth}
Let $\tG$, $\tH$ be graphs with non-trivial bold homology. Is the homology $\bH(\tG\times\tH)$ of the product non-trivial?
Is there a K\"unneth-like theorem for bold homology (and more generally for the \"uberhomology) with respect to some graph operation?
\end{q}
\begin{q}
The Euler characteristic of $\bH$ is the coefficient of $1$ in the (bi)graded Euler characteristic of the \"uberhomology; that is
\[ \ddot{\chi}_\tG(q,t) = \sum_{i,j,k\in \mathbb{N}} (-1)^{j}{\rm rank} (\UH_{i,k}^{j}(\tG))q^{i}t^{k} \in \bZ[q,t]\ .\]
Can we recover other known graph invariants from $\ddot{\chi}$? More generally, is the \"uberhomology a categorification of some known graph polynomial?
\end{q}

\begin{q} Can we find a graph $\tG$ such that $\chi(\bH(\tG)) = 0$ and $\bH(\tG) \neq 0$?
\end{q}

\bibliographystyle{alpha}
\bibliography{bibliography}

\end{document}